\documentclass[12pt]{article}

\usepackage{amsmath}
\usepackage{amssymb}
\usepackage{amsthm}
\usepackage{amsfonts}
\usepackage{a4}
\usepackage{graphicx}
\usepackage{ esint }

\newtheorem{theorem}{Theorem}[section]
\newtheorem{lemma}[theorem]{Lemma}

\newtheorem{assu}[theorem]{Assumption}
\newcommand{\Proof}{\par\noindent{\em Proof. }}
\newcommand{\eop}{\nopagebreak\hspace*{\fill}$\Box$\smallskip}

\theoremstyle{remark}
\newtheorem*{rem}{Remark}

\newcommand{\N}{\Bbb N}
\newcommand{\Z}{\Bbb Z}
\newcommand{\R}{\Bbb R}

\def\calL{\mathcal{L}}

\def\calZ{\mathcal{Z}}
\def\eps{\varepsilon}

\def\e{\mathbf{e}}

\def\dist{\operatorname{dist}}
\def\Xint#1{\mathchoice
   {\XXint\displaystyle\textstyle{#1}}%
   {\XXint\textstyle\scriptstyle{#1}}%
   {\XXint\scriptstyle\scriptscriptstyle{#1}}%
   {\XXint\scriptscriptstyle\scriptscriptstyle{#1}}%
   \!\int}
\def\XXint#1#2#3{{\setbox0=\hbox{$#1{#2#3}{\int}$}
     \vcenter{\hbox{$#2#3$}}\kern-.5\wd0}}

\def\dashint{\Xint-}

\usepackage{color}

\definecolor{dark_green}{rgb}{0, .5, 0}

\begin{document}

\begin{center}
\begin{Large}
{\bf {An analysis of crystal cleavage  
in the passage from atomistic models to continuum theory }}
\end{Large}
\end{center}

\begin{center}
\begin{large}
Manuel Friedrich\footnote{ Universit{\"a}t Augsburg, Institut f{\"u}r Mathematik, 
Universit{\"a}tsstr.\ 14, 86159 Augsburg, Germany. {\tt manuel.friedrich@math.uni-augsburg.de}}
and Bernd Schmidt\footnote{Universit{\"a}t Augsburg, Institut f{\"u}r Mathematik, 
Universit{\"a}tsstr.\ 14, 86159 Augsburg, Germany. {\tt bernd.schmidt@math.uni-augsburg.de}}\\
\end{large}
\end{center}

\begin{center}
\today
\end{center}
\bigskip

\begin{abstract}
We study the behavior of brittle atomistic models in general dimensions under uniaxial tension and investigate the system for critical fracture loads. We rigorously prove that in the discrete-to-continuum limit the minimal energy satisfies a particular cleavage law with quadratic response to small boundary displacements followed by a sharp constant cut-off beyond some critical value. Moreover, we show that the minimal energy is attained by homogeneous elastic configurations in the subcritical case and that beyond critical loading cleavage along specific crystallographic hyperplanes is energetically favorable. In particular, our results apply to mass spring models with full nearest and next-to-nearest pair interactions and provide the limiting minimal energy and minimal configurations. 
\end{abstract}
\bigskip

\begin{small}

\end{small}

\tableofcontents

\section{Introduction}

In spite of its importance in applications, a thorough understanding of the cleavage behavior of brittle crystals remains a challenging problem in theoretical mechanics. Subject to tensile boundary loads these crystals undergo an elastic deformation for very small boundary displacements after which they fracture, typically in the form of cleavage along crystallographic hyperplanes. In particular, there is (almost) no plastic deformation in between the elastic and fracture regimes. Ab initio calculations (see, e.g., \cite{JarvisHayesCarter:01,HayesOrtizCarter:04}) show that in fact irreversible cleavage occurs after two parts of a specimen have been separated by only a few angstroms. In engineering applications the consequences of such brittle fracture can be considerably more severe than those of ductile fracture which is characterized by an extensive plastic regime before a specimen finally breaks apart. In this context, it is noteworthy that many materials while ductile at high temperature become brittle below a critical temperature. Therefore, it is crucial to identify critical loads for failure and to analyze the shape of cracks that are formed beyond critical loading.

To tackle the problems arising in fracture mechanics with variational methods Francfort and Marigo \cite{Francfort-Marigo:1998} have introduced Griffith models leading to a minimization problem for energy functionals comprising elastic bulk terms in the intact regions of the body and surface contributions defined on a set of co-dimension one, called the \textit{jump set}. In such a continuum setting the observation that due to the crystalline structure certain crack geometries are preferred is modeled by anisotropic surface energies, see e.g.\ \cite{Alicandro-Focardi-Gelli:2000, Negri:2003}. It is a challenging problem and a field of active research to derive effective continuum models of this form from discrete systems. In particular, the analysis of atomistic models may lead to rigorous arguments proving that in the multidimensional framework it is energetically optimal to cleave the body along specific crystallographic hyperplanes. 

In the engineering literature such discrete systems had been analyzed computationally in \cite{JarvisHayesCarter:01,HayesOrtizCarter:04} and formally by renormalization group techniques in \cite{NguyenOrtiz:2002}. Braides, Lew and Ortiz \cite{Braides-Lew-Ortiz:06} then showed analytically that in the continuum limit the energy satisfies a certain cleavage law with a universal form independent of the specific choice of the interatomic potential. In all these models the crack geometry is pre-assigned and fracture may only occur along planes leading effectively to a one-dimensional problem. However, in order to understand the physical and geometrical cause for cleavage in the fracture regime it is indispensable to examine vectorial problems in more than one space dimension. 

In our previous work \cite{FriedrichSchmidt:2011} we analyzed a two dimensional model problem where the discrete system under consideration was given by the portion of a triangular lattice in a rectangular strip. Under uniaxial tension a cleavage law was derived with quadratic response for small boundary values followed by a sharp constant cut-off beyond some critical value (for further motivation and additional references we refer to this article.). The model was examined for validity and failure of crystal cleavage and  minimizers were classified being either homogeneous elastic deformations or configurations cleaved along an optimal crystallographic line in the subcritical and supercritical case, respectively. To the best of our knowledge this model is a first approach to higher dimensional problems being frame indifferent in the vector-valued arguments and fitting into the framework of \textit{free discontinuity problems} introduced by De Giorgi and Ambrosio \cite{DeGiorgi-Ambrosio:1988}, i.e.\ coming along without a priori assumptions on the crack geometry. Although having possible applications, e.g.\ in the modeling of the stability of brittle nanotubes, the limitation to two dimensions and a special reference configuration remains unsatisfactory. 

The aim of the present work is to extend the abovementioned results to arbitrary space dimensions and a general class of `cell energies' including well known mass-spring models where the pair interaction of neighboring atoms is modeled by potentials of Lennard-Jones type. It turns out that this analysis is considerably more involved and major difficulties have to be overcome, which were bypassed in \cite{FriedrichSchmidt:2011} due to the special assumptions of 1.\ a planar geometry, thus avoiding the possibility that cracks may concentrate on lower dimensional structures, 2.\ a triangular lattice of atoms, which led to an isotropic linearized elastic energy, and 3.\ nearest neighbor interactions that, up to symmetry, resulted in only one crack mode of a lattice triangle. Accordingly, in this case the basic energetic cleavage law (cf.\ \cite[Thm.~2.1]{FriedrichSchmidt:2011}) could be shown in a comparatively elementary way by resorting to slicing methods and convexity estimates in combination with a suitable projection technique. By way of contrast, the analogous result in arbitrary dimensions with general lattices and interaction potentials requires 1.\ new projection estimates for the size of cracks in the specimen, 2.\ a full dimensional analysis of an anisotropic mesoscopic auxiliary problem in various regimes and 3.\ a thorough analysis of all possible crack modes of a possible lattice unit cell.

To be more specific, in our model the atoms in the reference configuration are given by the portion $\eps\calL \cap \Omega$, where the macroscopic region $\Omega \subset \R^d$ occupied by the body is a cuboid and $\eps\calL$ is some Bravais lattice scaled by the typical interatomic distance $\eps \ll 1$. The main structural assumption is that the energy of a deformation $y: \eps\calL \cap \Omega \to \R^d$  may be decomposed as a sum over cell energies. The cell energy on its part depends on the discrete gradient $\bar{\nabla} y$ encoding all the relative displacements of atoms in a cell and satisfies some reasonable assumptions in the elastic regime (see e.g.\ \cite{ContiDolzmannKirchheimMueller:06}), particularly the frame indifference. As forces between well separated atoms are governed by dipole interactions we assume that for large deformation gradients the cell energy reduces to a pair interaction energy neglecting multiple point interactions. In the language of Truskinovsky's seminal contribution \cite{Truskinovsky:1996} our cell energy thoroughly describes the atomistic interactions within the two `phases' of brittle fracture: the linearly elastic and the fracture regime. 

Our main result is that under tensile boundary conditions the limiting minimal energy satisfies a cleavage law of a universal form essentially only depending on the stiffness and toughness of the material which may be deduced from the cell energy. In particular we confirm previous results in the mechanics and mathematics literature for one-dimensional models comprising effective interplanar potentials, cf.\ \cite{HayesOrtizCarter:04,NguyenOrtiz:2002,Braides-Lew-Ortiz:06}. Suitable test configurations then show that asymptotically optimal configurations are given by homogeneous elastic deformations showing the Poisson effect for subcritical boundary values and by configurations cleaved along specific crystallographic hyperplanes beyond critical loading. 

Deriving effective continuum theories for brittle materials is a challenging task and still an open problem in its full generality. One major complication is the possible coexistence and interaction of two different competing energy forms. In particular, the crack geometry might become extremely complex due to relaxation of the elastic energy by oscillating crack paths and infinite crack patterns occurring on different scales. 

In the present context of a uniaxial tension test which is a natural setting for the investigation of cleavage phenomena we overcome this difficulty by solving an auxiliary problem on a `mesoscopic cell' whose size is carefully chosen between the microscopic scale $\eps$ and the macroscopic magnitude of the specimen. On the one hand, by choosing this size small enough it is possible to separate the effects arising from the bulk elastic and the surface crack energy and to apply elaborated methods in the various regimes, including rigidity estimates \cite{FrieseckeJamesMueller:02} and slicing techniques (see e.g.\ \cite{Ambrosio-Fusco-Pallara:2000}).  On the other hand, given that the size is large with respect to $\eps$ we can exploit the validity of the Cauchy-Born-rule for sufficiently small strains which means, loosely speaking, that every single atom follows the mesoscopic deformation gradient and atomistic oscillations are effectively excluded (see \cite{ContiDolzmannKirchheimMueller:06, FrieseckeTheil:02}). More precisely and in mathematical terms, passing simultaneously from discrete to continuum theory and from finite to infinitesimal elasticity the discrete gradient of the atomic displacements reduces to a classical gradient leading to a simpler description of the stored elastic energy (cf.\ \cite{Schmidt:2009}). Moreover, with the help of tailor-made interpolations depending on individual crack modes, it can be shown that the fracture energy consisting of all contributions from pair interactions of neighboring atoms reduces to a surface energy in the continuum limit which only depends on the crack geometry (cf.\ also \cite{Braides-Gelli:2002-2}) and is minimized for a specific crystallographic hyperplane. Finally, we also show that intermediate regimes are energetically unfavorable. This is accomplished by establishing a $p$-growth estimate from below for carefully chosen $p > 1$. The analysis of the mesoscopic problem is furthermore complicated by the fact that individual cells on the boundary, where boundary values are prescribed, might have fractured. We resolve this problem by providing some estimates on the length of Lipschitz curves in sets of finite perimeter. 
       
We note that a complete characterization of minimizing sequences and proving strong discrete-to-continuum convergence results as carried out in the two-dimensional model problem \cite{FriedrichSchmidt:2011} seems currently out of reach. Moreover, the incorporation of more general boundary conditions would be desirable including the case of uniaxial compression which is as the uniaxial tension test a natural problem. It seems, however, that our techniques do not apply in the compressive case and that possibly additional modeling assumptions are necessary concerning an adequate and physically reasonable impenetrability condition. Another question which is beyond the scope of the present work is to devise models for ductile fracture on the microscopic scale that lend themselves to a similar analysis and that in the continuum limit will show an elastic as well as a plastic regime before ultimately leading to rupture. Nevertheless, we believe that the methods in this work, in particular the introduction of a mesoscopic scale whereby elastic and fracture effects can be separated, may contribute to solve more general problems in the future.

The paper is organized as follows. In Section 2 we introduce the discrete model and state the main cleavage law result. By a heuristic argument we determine the most interesting regime of boundary values, namely the one where the energies of typical elastic deformations and configurations with cleavage are of the same order. This scaling was first proposed by Nguyen and Ortiz \cite{NguyenOrtiz:2002}, who investigated the problem with renormalization group techniques.

Section 3 is devoted to preliminaries. We first derive formulae for the essential constants appearing in the cleavage law characterizing the stiffness and the toughness of the material. Here we already see that it is optimal to cleave along a crystallographic hyperplane. We introduce interpolations both for the elastic regime following the ideas in \cite{Schmidt:2009} and for the fracture regime being adapted for the application of slicing techniques. Moreover, we recall the definition and fundamental properties of special functions of bounded variation and state a short lemma about the length of Lipschitz curves in sets of bounded variation being substantially important in dimensions $d \ge 4$. 

Section 4 contains the essential technical estimates providing a lower comparison potential for the energy of a `cell of  mesoscopic size' under given averaged boundary conditions. The proof is mainly divided into three parts each of which dealing with one particular regime: The elastic regime where we show that linear elasticity theory applies, the fracture regime where we use a slicing argument in the framework of $SBV$ functions and an intermediate regime. Beyond that, in the case $d\ge 4$ an additional intermediate regime has to be introduced due to the fact that in higher dimensions it becomes more difficult to derive uniform bounds on the difference of boundary values. 

Section 5 is devoted to the proof of the main theorem which relies on the application of the comparison energy derived in Section 4 and a slicing argument in the space direction were the tensile boundary conditions were imposed. 

Finally, in Section 6 we give some examples of mass-spring models to which the aforementioned results apply and provide the limiting minimal energy as well as asymptotically optimal configurations. We first discuss the nearest neighbor interaction in a triangular lattice re-deriving results established in \cite{FriedrichSchmidt:2011} and also analyze the nearest and next-to-nearest neighbor interaction in a square lattice. In the latter case we see that in addition to the Poisson effect elastic minimizers generically also show a shear effect due to the anisotropy of the linearized elastic energy. Whereas the energetically favorable crack line in the triangular lattice was exclusively determined by the geometry of the problem, we find that for the square lattice two competing crystallographic lines occur due to possible different microscopic structures of fracture. Finally, we apply our results to a general nearest and next-to-nearest neighbor model in 3D considered e.g.\ in \cite{Schmidt:2006, Schmidt:2009}.

\section{The model and main results}\label{sec:model-and-main-results}

\subsection{The discrete model}\label{sec: model}

Let $\Omega \subset \R^d$ be the macroscopic region occupied by the body under consideration. To simplify the exposition we assume that $\Omega = (0,l_1) \times \ldots \times (0,l_d)$ is rectangular, but remark that all our results extend without difficulty to more general geometries as $\Omega = (0, l_1) \times \omega$, $\omega \subset \R^{d-1}$ open, for which cleavage boundary values as discussed in Section \ref{sec: bc} below may be imposed. Let ${\cal L}$ be some Bravais lattice in $\R^d$, i.e.\ there are linearly independent vectors $v_1, \ldots, v_d \in \R^d$ such that
$${\cal L} = \left\{\lambda_1 v_1 + \ldots \lambda_d v_d: \lambda_1,\ldots,\lambda_d \in \Z\right\} = A\Z^d,$$
where $A$ is the matrix $(v_1,\ldots,v_d)$. Without restriction we may assume that the vectors $v_i$ are labeled such that $\det A > 0$. The portion of the scaled lattice ${\cal L}_\eps = \eps {\cal L}$ lying in $\Omega$ represents the positions of the specimen's atoms in the reference position. Here $\eps$ is a small parameter measuring the typical interatomic distance eventually tending to zero. Note that  ${\cal L}_\eps$ partitions $\R^d$ into cells of the form $\eps A (\lambda + [0,1)^d)$ for $\lambda \in \Z^d$. The shifted lattice $\eps A ((\frac{1}{2},\ldots,\frac{1}{2})^T + \Z^d)$ consisting of the midpoints of the cells is denoted by $\calL'_\eps$. For $x \in \R^d$ we denote by $\bar{x} = \bar{x}(x,\eps)$ the center of the $\eps$-cell containing the point $x$ and set $Q_\eps(x) = \bar{x}(x,\eps) + \eps A [-\frac{1}{2}, \frac{1}{2})^d$. 

We choose a numbering $z_1, \ldots, z_{2^d}$ of the corners $A\left\{-\frac{1}{2},\frac{1}{2}\right\}^d$ of the reference cell $A[-\frac{1}{2},\frac{1}{2})^d$ and set 
\begin{align}\label{eq: Z def} 
Z = (z_1, \ldots, z_{2^d}), \ \ \ {\cal Z} = \left\{z_1,\ldots, z_{2^d}\right\}.
\end{align}
For subsets $U \subset \Omega$ we define the following lattice subsets with respect to the midpoints ${\cal L}'_\eps$ and the corners ${\cal L}_\eps$:
\begin{align*}
\calL'_\eps(U) =  \left\{\bar{x} \in \calL'_\eps: Q_\eps(\bar{x}) \cap U \neq \emptyset\right\}, \ \ 
\calL_\eps(U) = \calL'_\eps(U) + \eps\left\{z_1,\ldots,z_{2^d}\right\}, \\
(\calL'_\eps(U))^\circ =  \left\{\bar{x} \in \calL'_\eps: \overline{Q}_\eps(\bar{x}) \subset U\right\}, \ \ 
(\calL_\eps(U))^\circ = (\calL'_\eps(U))^\circ + \eps\left\{z_1,\ldots,z_{2^d}\right\}. 
\end{align*}
We call $Q_\eps(\bar{x})$ for $\bar{x} \in (\calL'_\eps(\Omega))^{\circ}$ an inner cell and set $\Omega_\eps = \bigcup_{\bar{x} \in (\calL'_\eps(\Omega))^{\circ}} Q_\eps(\bar{x})$. 

The deformations of our system are mappings $y :  \calL_\eps \cap \Omega \to \R^d$. Given $x \in \Omega_\eps$ and the corresponding midpoint $\bar{x} \in (\calL'_\eps(\Omega))^\circ$ we denote the images of the atoms in $\overline{Q_\eps(x)}$ by $y_i = y(\bar{x} + \eps z_i)$ for $i=1,\ldots, 2^d$ and view
\begin{align}\label{eq: y def}
Y(x) = (y_1, \ldots, y_{2^d}).
\end{align}
as elements of $\R^{d \times 2^d}$. We define the discrete gradient $\bar{\nabla} y(x) \in \R^{d \times 2^d}$ by
\begin{align}\label{eq: discrete grad}
\bar{\nabla} y := \eps^{-1} (y_1 - \bar{y}, \ldots, y_{2^d} - \bar{y}), \ \ \ \bar{y} := \frac{1}{2^d} \sum^{2^d}_{i=1} y_i.
\end{align}
In particular, $\bar{\nabla} y$ is a function on $\Omega_\eps$, which is constant on each cube $Q_\eps(\bar{x})$, $\bar{x} \in (\calL'_\eps(\Omega))^\circ$.

We also need to keep track of the atomic positions within subsets of ${\cal Z}$. Therefore, for a given matrix $G= (g_1, \ldots, g_{2^d}) \in \R^{d \times 2^d}$ and $\tilde{\calZ} \subset \calZ$ we define
\begin{align}\label{eq: G def}
G[\tilde{\calZ}] = \left(g_{j}\right)_{z_j \in \tilde{\calZ}} \in \R^{d \times \# \tilde{\calZ}}.
\end{align}
In cells with large deformation it will be convenient to measure the distance of different subsets of the atoms forming the cell.  For $G \in \R^{d \times 2^d}$ and ${\cal Z}_1, {\cal Z}_2 \subset \calZ$ we set
\begin{align}\label{eq: d def}
d(G; \calZ_1, \calZ_2) := \min\left\{|g_i - g_j|: z_i \in \calZ_1, z_j \in {\cal Z}_2\right\}.
\end{align}
We now define the set of interaction directions 
$$  {\cal V} = A \lbrace-1,0,1\rbrace^d \setminus \{0\}$$
and characterize the crystallographic hyperplanes spanned by the corners of a unit cell by their normal vectors. Let $S^{d-1} = \lbrace\xi \in \R^d: |\xi| = 1\rbrace$ and set 

\begin{align*}
{\cal P} &:= \lbrace\xi \in S^{d-1}: \exists u_1, \ldots, u_{d-1} \in {\cal V}, \ \text{span} \lbrace u_1, \ldots u_{d-1}\rbrace = \xi^{\perp} \rbrace.
\end{align*}
Note that every hyperplane is represented twice in ${\cal P}$, by $\xi$ and $-\xi$.
 
Our basic assumption is that the energy associated to deformations $y: \calL_\eps \cap \Omega \to \R^d$ can be written as a sum over cell energies $W_{\rm cell} : \R^{d \times 2^d} \to [0,\infty]$ in the form 
\begin{align}\label{eq: decomposition}
E_\eps (y) = \sum_{\bar{x} \in (\calL'_\eps(\Omega))^\circ} W_{\rm cell} (\bar{\nabla} y (\bar{x})).
\end{align}
For convenience the energy is defined as a sum over the inner cells only as the energy contribution of cells with midpoints lying in $\calL'_\eps(\Omega) \setminus (\calL'_\eps(\Omega))^\circ$ are negligible in our model for uniaxial extension. We briefly note that $W_{\rm cell}$ is of order one in atomic units and therefore we will have to consider a suitably scaled quantity of $E_\eps$ to arrive at macroscopic energy expressions for small $\eps$. This will be discussed in the next section.

\begin{rem} 
A decomposition as in \eqref{eq: decomposition} is, in particular, possible for many mass spring models, as will be exemplified in Section \ref{sec: examples}: The energy stored in an atomic bonds which lies on a face of more than one unit cell will then be equidistributed to the energy contribution of all adjacent cells. Moreover, energy functionals of the form \eqref{eq: decomposition} can also incorporate bond angle dependent energy terms. 
\end{rem}

We let 
$$\bar{SO}(d) := \left\{\bar{R} = RZ: R \in SO(d)\right\} \subset \R^{d  \times 2^d},$$
where $Z$ is as defined in \eqref{eq: Z def} and now describe the general assumptions on the cell energy $W_{\rm cell}$ in detail. 

\begin{assu}\label{assu: 1}
\begin{enumerate}
	\item[(i)] $W_{\rm cell}: \R^{d \times 2^d} \to [0,\infty]$ is invariant under translations and rotations, i.e.\ for $G \in \R^{d \times 2^d}$ we have
	      $$W_{\rm cell}(G) = W_{\rm cell} (RG + (c,\ldots,c)) $$
	    for all $R \in SO(d)$ and $c \in \R^d$.
	\item[(ii)] $W_{\rm cell}(G) = 0$ if and only if there exists $R \in SO(d)$ and $c \in \R^d$ such that
	      $$G = RZ + (c,\ldots,c). $$
	\item[(iii)] $W_{\rm cell}$ is continuous and $C^2$ in a neighborhood of $\bar{SO}(d)$. The Hessian $Q_{\rm cell} = D^2W_{\rm cell}(Z)$ at                  the identity is positive definite on the complement of the subspace spanned by translations $(c,\ldots,c)$ and                            infinitesimal rotations $HZ$, with $H + H^T = 0$.
	\item[(iv)] If there is a partition $\calZ = \calZ_1 \dot{\cup} \ldots \dot{\cup} \calZ_n$ such that $\min_{1\leq i < j \leq n} d(G;\calZ_i,\calZ_j)$ is near infinity for $G\in \R^{d \times2^d}$ then the energy $W_{\rm cell}$ decomposes, i.e.\ there are $W^{\calZ_i}: \R^{d \times \# \calZ_i} \to [0,\infty)$, $i=1,\ldots,n$, and $\beta(z_s,z_t) = \beta(z_t,z_s) \ge 0$ for $z_s,z_t \in \calZ$, $z_s \neq z_t$, such that 
	              $$W_{\rm cell}(G) = \sum^n_{i=1} W^{\calZ_i}(G[\calZ_i]) + \frac{1}{2}\sum_{1\leq i, j \leq n \atop i \ne j} \sum_{z_s \in \calZ_i} \sum_{z_t \in \calZ_j} \beta(z_s,z_t) + o(1)$$
	             as $\min_{1\leq i < j \leq n} d(G;\calZ_i,\calZ_j) \to \infty$, where the triple sum on the right hand side is strictly positive unless $n = 1$. 
	            The components of the energy satisfy 
	             \begin{align}\label{eq: compatible}
	             W^{\calZ_i}(H[\calZ_i]) \leq C W_{\rm cell}(H)
	             \end{align}
	             in a neighborhood of $\bar{SO}(d)$.
	\end{enumerate}
\end{assu}
Note that the above assumptions imply that the quadratic form $Q_{\rm cell}$ satisfies
\begin{align}\label{eq: Qcell prop}
Q_{\rm cell}(c,\ldots,c) = 0 , \ \ \ Q_{\rm cell} (HZ) = 0
\end{align}
for all $c \in \R^d$ and $H \in \R^{d \times d}$ with $H + H^T = 0$. Moreover, we obtain 
\begin{align}\label{eq: Wcell}
\liminf_{|G| \to \infty} W_{\rm cell}(G)> 0.
\end{align}
Partitioning the set of interaction directions as ${\cal V} = {\cal V}_1 \cup \ldots \cup {\cal V}_d$ with 
$${\cal V}_k = \Big\{ A \, t: t \in \lbrace-1,0,1\rbrace^d, \# t = k \Big\}$$
for $1 \leq k \leq d$, where $\# t$ denotes the number of non zero entries of the vector $t \in \lbrace-1,0,1\rbrace^d$, we may assume that for all $1\leq k \leq d-1$ and for all $\nu \in {\cal V}_k$ there are $\beta(\nu) \ge 0$ such that
\begin{align}\label{eq: beta nu}
\frac{1}{2} \beta(z_s,z_t) =  2^{k-d}\beta(\nu) \text{   for all } z_s,z_t \in \calZ: z_s  - z_t = \nu.
\end{align}
Indeed, a bond in $\nu$-direction, $\nu \in {\cal V}_k$ is shared by $2^{d-k}$, $1\leq k \leq d-1$ different cells. If these cells give different interaction energies $\beta^1, \ldots, \beta^{2^{d-k}}$ we would set $\beta(\nu) = \frac{1}{2}(\beta^1 + \ldots + \beta^{2^{d-k}})$ and would replace $\beta^i$ by $2^{k-d+1} \beta(\nu)$ without affecting the energy. The additional factor of $\frac{1}{2}$ takes account of the fact that every atomic bond is represented twice in ${\cal V}$. 

\begin{rem} 
Pair interaction potentials fulfilling assumption (iv) are, e.g., potentials of `Lennard-Jones-type'. They are characterized by the fact that for large distances the interaction energy is near a fixed positive value. We note that in the case of large deformation the energy reduces to a pair interaction energy neglecting multiple point interactions. This is meaningful as the forces between well separated atoms are governed by dipole interactions, while angle-dependent potentials are needed in order to describe the chemical binding effects for close atoms appropriately. Condition \eqref{eq: compatible} is a compatibility condition and ensures that surface terms cannot dominate bulk terms (cf.\ \cite{Schmidt:2009}).
\end{rem}

\subsection{Boundary values and scaling}\label{sec: bc}

We are interested in the behavior of the specimen under uniaxial extension, say in $\mathbf{e}_1$-direction. In particular, we would like to determine the critical value of the boundary displacement at which minimizers are no longer elastic deformations but form cracks. Moreover, we investigate if a separation of the body along specific crystallographic hyperplanes is indeed energetically most favorable.
In order to avoid geometric artefacts and complicated crack geometries, we will therefore assume that the specimen is `long enough' so that it is possible for the body to completely break apart along crystallographic hyperplanes not passing through the boundary parts $B_1 = \left\{x \in \overline{\Omega}: x_1=0\right\}$ and $B_2 = \left\{x \in \overline{\Omega}: x_1=l_1\right\}$. 

More precisely, we choose $l_1 \geq L=L(\sqrt{A^T A},W_{\rm cell}, l_2, \ldots, l_d)$ (see \eqref{eq: L def} below). In particular, we will see that $L$ may be chosen independently of the orientation of the lattice. We briefly note that, under additional symmetry conditions on the cell energy, $L$ is independent of $W_{\rm cell}$ and thus the minimum length only depends on the geometry of the problem. Such a symmetry condition is satisfied, e.g., if springs associated to the set ${\cal V}_k$ are of `the same type', i.e.\ give the same interaction energy for large expansion (see the examples in Section \ref{sec: examples}). We mention that one may overcome this technical difficulty alternatively by imposing periodic boundary conditions or considering infinite crystals.   

Due to the discreteness of the underlying atomic lattice the boundary conditions of uniaxial extension have to be imposed in atomistically small neighborhoods of $B_1$ and $B_2$ as otherwise unphysical boundary effects may occur, in particular cracks near the boundary might become energetically more favorable. Define
\begin{align}\label{eq: lA}
l_A = \sum^d_{j=1} |v_j \cdot \e_1|
\end{align}
and let $B^\eps_1 = \left\{x \in \overline{\Omega}: x_1 \le 2 l_A\eps\right\}$ and $B^\eps_2 = \left\{x \in \overline{\Omega}: x_1 \geq l_1 - 2 l_A\eps\right\}$. For $a_{\eps} > 0$ we set 
\begin{align} \label{eq: bc1}
\begin{split}
  {\cal A}(a_{\eps}) 
  &= \big\{ y = (y^1,\ldots, y^d)^T : {\cal L}_\eps \cap \Omega \to \R^d : \\ 
  &\qquad \qquad \qquad \qquad  
     y^1(x) = (1 + a_{\eps}) x_1 \text{ for } x \in B^\eps_1 \cup B^\eps_2\big\}. 
     \end{split}
\end{align}
There is some arbitrariness in the implementation of boundary conditions. A possible alternative is, e.g.,
\begin{align}\label{eq: bc2}
  y^1(x) = x_1 \text{ for } x \in B^\eps_1 ~\text{ and }~ 
  y^1(x) = x_1 + a_{\eps} l_1 \text{ for } x \in B^\eps_2. 
\end{align}
We remark that such different choices do not change the results of our analysis and we will say that a deformation satisfies the boundary condition if either \eqref{eq: bc1} or \eqref{eq: bc2} is satisfied. Note that there are no assumptions on the other $d-1$ components of the deformation $y$ near the boundaries $B_1$ and $B_2$ of the boundary displacement, i.e.\ the atoms may `slide along the boundary'.  

There are two obvious choices for deformations satisfying the boundary conditions: The homogeneous elastic deformation $y^{\rm el}(x) =  (1+ a_{\eps}) x$ and a cracked body deformation $y^{\rm cr}$, which has the form $y^{\rm cr} = x \chi_{\Omega_1} + (x + a_\eps l_1 \e_1) \chi_{\Omega_2}$, where the sets $\Omega_1$ and $\Omega_2$ form a partition of $\Omega$ and are separated by some hyperplane (or manifold) intersecting $\overline{\Omega}$ in the set $\overline{\Omega} \setminus (B^\eps_1 \cup B^\eps_2)$. Noting that $W_{\rm cell} \sim \dist^2(\cdot,\bar{SO}(d))$ in a neighborhood of $\bar{SO}(d)$ and employing \eqref{eq: compatible} it is not hard to see that for $\eps \ll a_\eps \ll 1$ we have 
$$ E_{\eps}(y^{\rm el}) \sim \eps^{-d} a_\eps^2, \quad \quad 
   E_{\eps}(y^{\rm cr}) \sim \eps^{1-d}. $$ 
We are particularly interested in the regime where both of these energy values are of the same order, i.e.\ $E_{\eps}(y^{\rm el}) \sim E_{\eps}(y^{\rm cr}) \sim \eps^{1-d}$. This implies $a_\eps \sim \sqrt{\eps}$. As alluded to above in order to arrive at finite and nontrivial energies in the limit $\eps \to 0$, we rescale $E_{\eps}$ to ${\cal E}_{\eps} := \eps^{d-1} E_{\eps}$.

\subsection{Cleavage law}

We now state our main result about the limiting minimal energy as $\eps \to 0$ when $a_\eps/\sqrt{\eps} \to a \in [0,\infty]$. In particular, we will see that the minimal energy is given by elastic deformations for $a_\eps$ up to some critical value $a_{\rm crit}$ of boundary displacements and by cleavage along a specific crystallographic hyperplane beyond this value. 
Before we state the theorem we introduce two constants occurring in the limiting minimal energy which describe the stiffness and toughness of the material. The constant in the fracture energy is given by
\begin{align}\label{eq: beta def}
\beta_A = \min_{\xi \in \R^d {\setminus \{0\}}} \frac{\sum_{\nu \in {\cal V}} \beta(\nu) |\nu \cdot \xi|}{|\e_1 \cdot \xi|}
\end{align}
with $\beta(\nu)$ as in \eqref{eq: beta nu}. In Lemma \ref{lemma: betaA min} below we will show that the minimum $\beta_A$ is attained for some $\xi \in {\cal P}$. In particular, this means that cleavage along a crystallographic hyperplane is energetically favorable.  Concerning the elastic regime we define a \textit{reduced energy} for the quadratic form $Q_{\rm cell}$ by
\begin{align}\label{eq: Wred def}
\tilde{Q}(r) := \min\left\{Q_{\rm cell}(e(G)\cdot Z): G \in \R^{d\times d}, g_{11} = r\right\}
\end{align} 
for $r \in \R$, $e(G) = \frac{1}{2}(G + G^T)$. As the problem is quadratic with a linear constraint, it is not hard so see that
$\tilde{Q}(r) = \alpha_A r^2 $ for a specific $\alpha_A > 0$. We will see that
\begin{align}\label{eq: alpha lax}
\alpha_ A =\frac{\det ({\cal Q})}{\det (\hat{\cal Q})},
\end{align}
where, roughly speaking, $\cal Q$ is the projection of $Q_{\rm cell}$ onto the linear subspace orthogonal to infinitesimal rotations and $\hat{\cal Q}$ arises from $\cal Q$ by cancellation of the first row and column. This will be stated more precisely in Lemma \ref{lemma: alpha Q} below.
\begin{theorem}\label{th: cleavage}
Let $l_1\geq L(\sqrt{A^T A},W_{\rm cell},l_2,\ldots,l_d)$ and suppose $a_\eps / \sqrt{\eps} \to a \in [0, \infty]$. The limiting minimal energy is given by
$$ {\cal E}_{\rm lim}(a) 
   :=  \lim_{\eps \to 0} \inf \lbrace {\cal E}_\eps (y)  : y \in {\cal A}(a_\eps)\rbrace =  \frac{\prod^d_{j=2} l_j}{\det A} \min \Big\lbrace \frac{1}{2} l_1 \alpha_A a^2, \beta_A \Big\rbrace.$$
\end{theorem} 
As discussed above, for $a \in \lbrace 0, \infty \rbrace$ either the elastic or the fracture regime is energetically favorable. The more interesting case is $a \in (0,\infty)$ where both energies are of the same order. The limiting minimal energy satisfies
a cleavage law of a universal form essentially only depending on the stiffness and toughness of the material and exhibits quadratic response to small boundary displacements followed
by a sharp constant cut-off beyond the critical value of boundary displacements 
\begin{align}\label{eq: a crit} 
a_{\rm crit} = \sqrt{\frac{2 \beta_A}{l_1 \alpha_A}}. 
\end{align}

We briefly indicate asymptotically optimal configurations. In the subcritical case $a \le a_{\rm crit}$ we consider the sequence of configurations
\begin{align*}
y^{\rm el}_\eps(x) = x +  \bar{F}(a_\eps) \, x, \ \ \ x \in \calL_{\eps} \cap \Omega,
\end{align*}
where $\bar{F}(a_\eps)$ is the solution of the minimization problem \eqref{eq: Wred def} with $r = a_\eps$ (see Lemma \ref{lemma: alpha Q} below). The deformations behave purely elastically and as we will see in the examples in Section \ref{sec: examples} show elongation in $\e_1$-direction  and contraction in the other space directions, a manifestation of the Poisson effect. Moreover, the configurations illustrate the validity of the Cauchy-Born-rule in this regime as each individual atom follows the macroscopic deformation gradient. In the supercritical case $a\ge a_{\rm crit}$ there is some $\xi \in {\cal P}$ and $c \in \R$ such that the hyperplane $\Pi = \lbrace x \in \R^d: x \cdot \xi = c\rbrace$ satisfies $\Pi \cap \overline{\Omega} \subset \overline{\Omega}\setminus (B_1^\eps \cup B^\eps_2)$ and the configurations
\begin{align}\label{eq: fracture min2}
  y^{\rm cr}_\eps(x) 
  = \begin{cases}  
       x, &  x\cdot \xi < c, \\
       x + l_1 a_\eps  \e_1, & x\cdot \xi > c, 
    \end{cases} \ \ \ x \in \calL_{\eps} \cap \Omega, 
\end{align}
are asymptotically optimal. As $\xi \in {\cal P}$, we conclude that $\Pi$ is a crystallographic hyperplane, as desired.  

Let us also remark that our class of atomistic interactions is rich enough to model any non-degenerate linearly elastic energy density, respectively, any preferred cleavage normal in ${\cal P}$ (not perpendicular to $\e_1$) in the continuum limit: In the elastic regime, this has in fact been observed in \cite[Prop.~1.10]{Braun-Schmidt:2013}. Now suppose that $\xi \in {\cal P}$ with $\xi \cdot \e_1 \ne 0$ is orthogonal to $\text{span} \{u_1, \ldots, u_{d-1}\}$ with $u_1, \ldots, u_{d-1} \in {\cal V}$. Then, if $\beta(u_1), \ldots, \beta(u_{d-1})$ are much larger than $\beta(\nu)$ for all $\nu \in {\cal V} \setminus \{u_1, \ldots, u_{d-1}\}$, it is elementary to see that the minimum of 
$$ \min_{\varsigma \in {\cal P}} \frac{\sum_{\nu \in {\cal V}} \beta(\nu) |\nu \cdot \varsigma|}{|\e_1 \cdot \varsigma|}, $$ 
is attained at $\varsigma = \xi$, so that indeed $\xi$  defines a crack normal for an asymptotically optimal configuration as in \eqref{eq: fracture min2}. 

The main idea in the proof of Theorem \ref{th: cleavage} is based on the derivation of a lower comparison potential for a certain cell energy depending on the expansion in $\e_1$-direction and on the application of a slicing argument. Testing either with elastic deformations or configurations forming jumps along specific hyperplanes as given above, we will then see that this lower bound is sharp. Actually, it will turn out that the lower bound coincides with the reduced energy $\tilde{Q}$ in the regime of infinitesimal elasticity.     
In contrast to the local definition of $\tilde{Q}$, however, it is in general not convenient to optimize the energy $W_{\rm cell}$ of single cells individually as it is geometrically nonlinear and therefore, due to possible rotations, the corresponding minimizer for one cell might not be compatible with deformations defined on the whole domain. As a remedy we will introduce a mesoscopic localization technique and will consider `large cells' defined on a mesoscopic scale $\eps^{\frac{3d-1}{3d}}$. This main technical result is addressed in Section \ref{sec: mesoscopic cell}.

\section{Preliminaries}\label{sec: pre}

\subsection{Elementary properties of the cell energy}

\subsubsection*{Elastic energy}

We first provide a lower bound for the cell energy. For that purpose, let $V_0 = \R^d \otimes (1, \ldots 1)$ denote the subspace of infinitesimal translations $(x_1, \ldots, x_{2^d}) \to (v,\ldots,v)$, $v \in \R^d$. 

\begin{lemma}\label{lemma: energy well}
For every $T > 0$ there is a constant $C>0$ such that
$$\dist^2(G,\bar{SO}(d)) \leq CW_{\rm cell}(G)$$
for all $G \in \R^{d \times 2^d}$ with $G \perp V_0$ and $|G| \leq T$.
\end{lemma}

\Proof The proof is essentially contained in  \cite[Lemma~3.2]{Schmidt:2006} and relies on the growth assumptions on $W_{\rm cell}$ near $\bar{SO}(d)$ (see Assumption \ref{assu: 1}(iii)). \eop

We now give a precise characterization of $\alpha_A$ (see \eqref{eq: alpha lax}). We may view a symmetric matrix $F = (f_{ij}) \in \R^{d \times d}_{\rm sym}$ as a vector $f = (f_1, \ldots, f_{\hat{d}}) \in \R^{\hat{d}}$, $\hat{d} = \frac{d(d+1)}{2}$, whose components are the entries $f_{ij}$ with $i \le j$, numbered such that $f_1 = f_{11}$. Then 
$$ Q_{\rm cell}(F \cdot Z) = f^T {\cal Q} f $$ 
for some symmetric positive definite ${\cal Q} \in \R^{\hat{d} \times \hat{d}}_{\rm sym}$. For each $r \in \R$ there is consequently a unique $\bar{f} \in \R^{\hat{d}}$ minimizing 
$$ f^T {\cal Q} f \quad \text{subject to} \quad f_1 = r $$ 
and a corresponding Lagrange multiplier $\lambda \in \R$ such that $\bar{f}^T {\cal Q} = \lambda \mathbf{e}_1^T$, i.e.\  
$$ \bar{f} = \lambda {\cal Q}^{-1} \mathbf{e}_1, $$
where $\mathbf{e}_1$ denotes the first canonical unit vector in $\R^{\hat{d}}$. Multiplying with $\mathbf{e}_1^T$ and using $\bar{f}_1 = r$, we find that $\lambda = \frac{r}{\mathbf{e}_1^T {\cal Q}^{-1} \mathbf{e}_1}$ and thus 
$$ \bar{f} = \frac{r}{\mathbf{e}_1^T {\cal Q}^{-1} \mathbf{e}_1} {\cal Q}^{-1} \mathbf{e}_1. $$
For the minimal value we obtain
$$ \tilde{Q}(r) 
   = \bar{f}^T {\cal Q} \bar{f} 
   = \frac{r^2}{\mathbf{e}_1^T {\cal Q}^{-1} \mathbf{e}_1}. $$
With $\hat{\cal Q}$ denoting the $(\hat{d} - 1) \times (\hat{d} - 1)$ matrix obtained form ${\cal Q}$ by deleting the first row and the first column and using that ${\cal Q}^{-1} = \frac{1}{\det {\cal Q}} \operatorname{cof} {\cal Q}$ and thus ${\cal Q}^{-1} \mathbf{e}_1 = \frac{1}{\det {\cal Q}} (\operatorname{cof} {\cal Q})_{\cdot 1}$ and $\mathbf{e}_1^T {\cal Q}^{-1} \mathbf{e}_1 = \frac{(\operatorname{cof} {\cal Q})_{11}}{\det {\cal Q}} = \frac{\det \hat{\cal Q}}{\det {\cal Q}}$, this can alternatively be written as 
$$ \bar{f} 
   = \frac{r \det {\cal Q}}{\det \hat{\cal Q}} (\operatorname{cof} {\cal Q})_{\cdot 1}, \qquad 
   \tilde{Q}(r) 
   = \frac{r^2 \det {\cal Q}}{\det \hat{\cal Q}}. $$

We summarize these observations in the following lemma:

\begin{lemma}\label{lemma: alpha Q}
The reduced energy satisfies 
$$ \tilde{Q}(r) = \alpha_A r^2 
   \quad \text{with} \quad 
   \alpha_A = \frac{1}{\mathbf{e}_1^T {\cal Q}^{-1} \mathbf{e}_1} = \frac{\det {\cal Q}}{\det \hat{\cal Q}}. $$ 
For each $r$ there exists a unique $\bar{F}(r) \in \R^{d \times d}_{\rm sym}$ which satisfies $\tilde{Q}(r) = Q_{\rm cell}(\bar{F}(r) \cdot Z)$ and $f_{11} = r$. $\bar{F}(r)$ depends linearly on $r$. 
\end{lemma}

\subsubsection*{Fracture energy}

In Theorem \ref{th: cleavage} we have seen that the limiting minimal fracture energy has the form \eqref{eq: beta def}. We now investigate this term in detail and determine the minimizers. We let 
\begin{align}\label{eq: crack}
\Lambda(\varsigma) := \frac{\sum_{\nu \in {\cal V} } \beta(\nu) |\nu \cdot \varsigma | - |\e_1 \cdot \varsigma|\beta_A}{|\varsigma|} 
\end{align}
and observe $\Lambda(\varsigma) \geq 0$ for all $\varsigma \in \R^d \setminus \{0\}$. Note that the minimum in the definition \eqref{eq: beta def} of $\beta_A$ and the minumum of $\Lambda$ in \eqref{eq: crack} are attained on the compact set $S^{d-1} = \lbrace\varsigma \in \R^d: |\varsigma| = 1\rbrace$. Moreover, we note that the minimizers in \eqref{eq: beta def} are precisely the minimizers of $\Lambda$ (with minimal value $0$). They are obviously not perpendicular to $\e_1$. 

\begin{lemma}\label{lemma: betaA min}
The minimum ($=0$) of $\Lambda(\xi)$ in $S^{d-1}$ is attained for some $\xi \in {\cal P}$. 
\end{lemma}

\begin{rem}
A natural guess would be that in fact every minimizer of $\Lambda$ lies in ${\cal P}$. Surprisingly this turns out to be wrong in general. There are (non-generic) models even leading to a continuum of optimal crack directions. This will be illustrated in Section \ref{sec: examples} for a basic mass spring model in 2d. As a consequence, for such a model it is not possible to prove that in the fracture regime the body has to break apart along crystallographic hyperplanes. 
\end{rem} 

\Proof For $\delta > 0$ we define 
$$\Lambda_\delta(\varsigma) = \Lambda(\varsigma) + \delta \frac{ |\varsigma \cdot\e_1|}{|\varsigma|}.$$
Obviously $\Lambda_\delta$ attains its minimum with $0 < \min_{\xi \in \R^d} \Lambda_\delta(\xi) \leq \delta$. We show that for small $\delta$ if a $\varphi \in S^{d-1}$ satisfies $\Lambda_\delta(\varphi)= \min_{\xi \in \R^d} \Lambda_\delta(\xi)$ then $\varphi  \in {\cal P}$. 

We define ${\cal U} = \lbrace\nu \in {\cal V}: \nu \cdot \varphi > 0\rbrace \cup \lbrace\e_1\rbrace$ and ${\cal U}_0 = \lbrace\nu \in {\cal V}: \nu \cdot \varphi = 0\rbrace$. Note that by \eqref{eq: crack} and \eqref{eq: beta def} $\varphi \cdot \e_1 \ne 0$, so without loss of generality we may assume that $\varphi \cdot \e_1 > 0$. For $\nu \in {\cal V} \setminus \lbrace\e_1\rbrace$ let $\tilde{\beta}(\nu) = 2 \beta(\nu)$. If $\e_1 \in {\cal V}$, we set $\tilde{\beta}(\e_1) =  2 \beta(\e_1) - \beta_A + \delta$. Otherwise we only set $\tilde{\beta}(\e_1) = - \beta_A + \delta$. If the claim were false, then $\text{dim}\,\text{span}\ {\cal U}_0 < d-1$. Therefore, we can choose some $\eta \in \R^d \setminus \{0\}$ such that $\eta \cdot \varphi = 0$ and $\eta \cdot \nu = 0$ for all $\nu \in {\cal U}_0$. We now investigate the behavior of $\Lambda_\delta$ at $\varphi$ in direction $\eta$. Using that $\nu \cdot \varphi = \nu \cdot \eta = 0$ for all $\nu \in {\cal U}_0$, for $|t|$ sufficiently small we obtain 
$$ \lambda(t) 
   := \Lambda_\delta(\varphi + t\eta) 
   = \frac{\sum_{\nu \in {\cal U}} \tilde{\beta}(\nu) \nu \cdot (\varphi + t\eta)}{|\varphi + t \eta|}.$$ 
We differentiate and obtain from $\eta \cdot \varphi = 0$
\begin{align*}
  \lambda'(t) 
  &= \frac{\sum_{\nu \in {\cal U}} \tilde{\beta}(\nu) \nu \cdot \eta}{|\varphi + t\eta|} 
   - \frac{\big(\sum_{\nu \in {\cal U}} \tilde{\beta}(\nu) \nu \cdot (\varphi + t\eta)\big)(\varphi + t\eta) \cdot \eta}{|\varphi + t\eta|^3} \\
  &=  \frac{\sum_{\nu \in {\cal U}}\tilde{\beta}(\nu) \nu \cdot \eta}{|\varphi + t\eta|} -  t |\eta|^2 \frac{\sum_{\nu \in {\cal U}} \tilde{\beta}(\nu)\nu \cdot (\varphi + t\eta)}{|\varphi + t\eta|^3}.
\end{align*}
If $\lambda'(0) \neq 0$ then $\varphi$ is not a critical point of $\Lambda_\delta$ which contradicts the above assumption. So we may assume that $\lambda'(0) = \frac{\sum_{\nu \in {\cal U}}\beta(\nu) \nu \cdot \eta}{|\varphi|} = 0$ and thus 
\begin{align*}
  \lambda'(t) 
  &= - t |\eta|^2 \frac{\sum_{\nu \in {\cal U}} \tilde{\beta}(\nu)\nu \cdot \varphi}{|\varphi + t\eta|^3}, 
\end{align*}
leading to the contradiction 
\begin{align*}
  \lambda''(0) = - |\eta|^2 \Lambda_{\delta}(\varphi) < 0. 
\end{align*}

Thus, we have shown that $\Lambda_\delta$ attains its minimum for some $\varphi \in {\cal P}$ for all $\delta > 0$. As $\min_{\xi \in \R^d} \Lambda_\delta(\xi) \leq \delta$,  passing to the limit $\delta \to 0$ we obtain the claim.
\eop

We are now in a position to render more precisely the definition of the minimum length assumed in Section \ref{sec: bc}. Let
$$M_1 = \min_{\xi \in S^{d-1}} \sum_{\nu \in {\cal V}} \beta(\nu) |\nu \cdot \xi|, \ \  M_2 = \max_{\xi \in S^{d-1}} \sum_{\nu \in {\cal V}} \beta(\nu) |\nu \cdot \xi|.$$
It is not hard to see that $M_1,M_2$ are independent of the particular rotation of the lattice. Then $\beta_A \leq M_2$ and therefore the minimizer $\xi \in S^{d-1}$ of \eqref{eq: beta def} satisfies $|\xi \cdot \e_1| \geq \frac{M_1}{M_2}$. Consequently, an elementary argument shows that choosing $C>0$ large enough independently of $M_1,M_2$, $l_2, \ldots, l_d$ and setting
\begin{align}\label{eq: L def}
L= L( \sqrt{A^T A},W_{\rm cell}, l_2,\ldots, l_d) = C \max\lbrace l_2,\ldots,l_d\rbrace \, \frac{M_2}{M_1}
\end{align}
we find that for specimens with $l_1 > L$ it is possible to completely break apart along hyperplanes not passing through the boundary parts $B_1$ and $B_2$.

\subsection{Interpolation}\label{sec: interpol}

In the following it will be useful to choose a particular interpolation $\tilde{y}$ of the lattice deformation $y:  \calL_\eps \cap U \to \R^d$ for $U \subset \Omega$ open. We introduce a threshold value $C_{\rm int} \geq 1$ to be specified later and first consider some cell $Q_\eps(\bar{x})$, $\bar{x} \in (\calL'_\eps(U))^\circ $,  where the lattice deformation satisfies $\text{diam}\left\{\bar{\nabla} y (\bar{x}), {\cal Z}\right\} \leq C_{\rm int}$. Here for $G = (g_1, \ldots, g_{2^d}) \in \R^{d \times 2^d}$ and $\tilde{\cal Z} \subset {\cal Z}$ we define $\text{diam}\left\{G, \tilde{\cal Z}\right\} := \max\left\{ |g_i - g_j|: z_i, z_j \in \tilde{\cal Z}\right\}$, so particularly we have 
\begin{align}\label{eq: max def}
\text{diam}\left\{\bar{\nabla} y (\bar{x}), \tilde{\cal Z}\right\} = \frac{1}{\eps}\max\left\{ |y(\bar{x} + \eps z_i) - y(\bar{x} + \eps z_j)|: z_i, z_j \in \tilde{\cal Z}\right\},
\end{align}
where $\bar{\nabla} y (\bar{x})$ is given in \eqref{eq: discrete grad}. We will call cells with this property `intact cells' and by ${\cal C}'_\eps \subset (\calL'_\eps(U))^\circ  $ we denote the set of their midpoints. The complement $\bar{\cal C}'_\eps := (\calL'_\eps(U))^\circ  \setminus {\cal C}'_\eps$ labels the centers of cells we consider to be `broken'. 

We first consider $Q_\eps(\bar{x})$ for $\bar{x} \in {\cal C}'_\eps$. The interpolation we use was introduced in \cite{Schmidt:2009}. We repeat the procedure here for the sake of completeness. Consider the reference cell $Q = A[-\frac{1}{2},\frac{1}{2})^d$ with deformation $y: A\left\{-\frac{1}{2},\frac{1}{2}\right\}^d = {\cal Z} \to \R^d$. We first interpolate linearly on the one-dimensional faces of $\overline{Q}$, which are given by segments $[z_i,z_j]$, where $z_i - z_j$ is parallel to one of the lattice vectors $v_n$, $n=1,\ldots,d$. Subsequently we consider two-dimensional faces and define a triangulation and interpolation as follows: Given a face $\text{co}\left\{z_{i_1}, z_{i_2}, z_{i_3}, z_{i_4}\right\}$ ($\text{co}$ denoting the convex hull) with 
$$z_{i_2} = z_{i_1} + v_n, \ \ \  z_{i_3} = z_{i_1} + v_n + v_m, \ \ \ z_{i_4} = z_{i_1} + v_m$$
we define 
$$\zeta = \frac{1}{4}(z_{i_1} + \ldots + z_{i_4}), \ \ \  y(\zeta) = \frac{1}{4}(y(z_{i_1}) + \ldots + y(z_{i_4}))$$ 
and interpolate linearly on each of the four triangles  $\text{co}\left\{z_{i_j},z_{i_{j+1}},\zeta\right\}$ for $j=1,\ldots,4$ with the convention $i_5=i_1$. In general, having chosen a simplicical decomposition as well as corresponding linear interpolations on the faces of dimension $n-1$ we decompose and interpolate on an $n$-dimensional face $F = \text{co}\left\{z_{i_1},\ldots, z_{i_{2^ n}}\right\}$ in the following way: Let
$$\zeta = \frac{1}{2^n} \sum^{2^n}_{j=1} z_{i_j}, \ \ \ y(\zeta) = \frac{1}{2^n} \sum^{2^n}_{j=1} y(z_{i_j}).$$
We decompose $F$ by the simplices $\text{co}\left\{w_1,\ldots,w_{n},\zeta\right\}$, where $\text{co}\left\{w_1,\ldots,w_n\right\}$ is a simplex belonging to the decompostion of an $(n-1)$-dimensional face constructed in a previous step. We now interpolate linearly on these simplices. 

For cells lying at the boundary $B^\eps_1$, $B^\eps_2$ we can repeat the above construction at least for the first component $y^1$: Let $\bar{x} \in \calL'_\eps(\Omega) \setminus (\calL'_\eps(\Omega))^\circ$ such that $\bar{x} + \eps z_i \notin \Omega$ implies $-l_A \eps \leq(\bar{x} + \eps z_i)_1\leq 0$ or $l_1\leq(\bar{x} + \eps z_i)_1\leq l_1 + l_A \eps$, respectively. Now let $y^1(\bar{x} + \eps z_i) = (1+ a_\eps) (\bar{x} + \eps z_i)_1$. Thus, the first component of $Y(\bar{x}) = (y_1,\ldots,y_{2^d})$ is well defined and we may proceed as above to construct $\tilde{y}^1$. 

We now concern ourselves with `broken cells' $\bar{x} \in \bar{\cal C}'_\eps$. Let $y$ be a corresponding lattice deformation defined on ${\cal Z}$. Recalling definition \eqref{eq: max def} we first choose a partition ${\cal Z}_1 \dot{\cup} \ldots \dot{\cup} {\cal Z} = {\cal Z}$ of the corners with 
$$\text{diam}\left\{\bar{\nabla} y (\bar{x}),{\cal Z}_i\right\} \leq \Big(\frac{\# {\cal Z}_i}{2^d}\Big)^2  C_{\rm int}.$$
We note that this partition can be chosen in a way that $\eps d(\bar{\nabla} y (\bar{x}); \calZ_i,\calZ_j) > 2^{-2d} \eps C_{\rm int}$ for all sets ${\cal Z}_i, \calZ_j$, $i \ne j$. Indeed, if there were $\bar{z}_i \in \calZ_i$, $\bar{z}_j \in \calZ_j$ such that $|y(\bar{z}_i) - y(\bar{z}_j)| = \eps d(\bar{\nabla} y (\bar{x}); \calZ_i,\calZ_j) \leq 2^{-2d} \eps C_{\rm int}$ then 
\begin{align*}
|y(z_i) - y(z_j)| 
&\leq |y(z_i) - y(\bar{z}_i)| + |y(\bar{z}_i) - y(\bar{z}_j)|  +|y(z_j) - y(\bar{z}_j)| \\
& \leq \frac{(\# \calZ_i)^2 + (\# \calZ_j)^2 + 1}{2^{2d}} \, \eps C_{\rm int} 
\leq \Big(\frac{\# (\calZ_i \cup \calZ_j)}{2^d}\Big)^2 \eps C_{\rm int}
\end{align*}
for all $z_i \in \calZ_i$, $z_j \in \calZ_j$  and we could set $\tilde{\calZ} = \calZ_i \cup \calZ_j$. Clearly, the cardinality of this partition is at least two for every cell $Q_\eps(\bar{x})$, $\bar{x} \in \bar{\cal C}'_\eps$. Then, by Assumption \ref{assu: 1}, particularly by \eqref{eq: Wcell}, it is not hard to see that there is a constant $C=C(C_{\rm int})$ such that 
\begin{align}\label{eq: least energy}
W_{\rm cell}(\bar{\nabla} y (\bar{x})) \geq C
\end{align}
for $\bar{x} \in \bar{\cal C}'_\eps$. Note that $C = C(C_{\rm int})$ can be chosen independently of $C_{\rm int}$ for $C_{\rm int} \ge 1$. 

We now choose the interpolation for each component $i$ separately as follows. If $\text{diam}_i\left\{\bar{\nabla} y(\bar{x}), {\cal Z}\right\} \leq C_{\rm int}$ we define $\tilde{y}^i$ as before for `intact cells' in ${\cal C}'_\eps \setminus \bar{\cal C}'_\eps$. Here for $G = (g_1, \ldots, g_{2^d}) \in \R^{d \times 2^d}$ we define similarly as in \eqref{eq: max def}
$$\text{diam}_i\left\{G, \cal Z\right\} := \max\left\{ |(g_j - g_k) \cdot \e_i|: z_j, z_k \in {\cal Z}\right\}.$$ Otherwise we set 
$$\tilde{y}^i(x) = y^i(z_1) + (x - z_1)_i $$
on $Q_\eps(\bar{x})$ and therefore $\nabla \tilde{y}^i = \e^T_i$. Consequently, $|\nabla \tilde{y}| \le CC_{\rm int}$ a.e.\  also on broken cells $Q_\eps(\bar{x})$, $\bar{x} \in \bar{\cal C}'_\eps$. Finally having constructed the interpolation on all cells $Q_\eps (\bar{x})$ we briefly note that $\tilde{y} \in SBV(\Omega_\eps, \R^d)$, cf.\ Section \ref{sec: sbv}. 

For every interaction direction $\nu \in {\cal V}$ we introduce a further interpolation $\bar{y}_\nu$ as follows. We choose vectors $\lbrace v_{i_1}, \ldots, v_{i_{d-1}}\rbrace \subset \lbrace v_1, \ldots, v_d\rbrace$ such that $\nu, v_{i_1}, \ldots, v_{i_{d-1}}$ are linearly independent and define the lattice ${\cal G}^\nu_\eps = \eps D^\nu \Z^d$ with $D^\nu = (\nu, v_{i_1}, \ldots, v_{i_{d-1}})$ partitioning $\R^d$ into cells of the form $Q^\nu_\eps(\lambda) := \eps D^\nu (\lambda + [0,1)^d)$ for $\lambda \in \Z^d$. Note that ${\cal L}_\eps = {\cal G}^\nu_\eps$. 

We describe the interpolation on the reference cell $Q^\nu = D^\nu [0,1)^d$ with deformation $y : D^\nu \lbrace0,1\rbrace^d \to \R^d$. If $|y(\nu) - y(0)| \leq  C_{\rm int}  \eps$ we let 
\begin{align}\label{eq: bar int def}
\bar{y}_\nu (D^\nu  x) = (1 - x_1) y(0) + x_1 y(\nu)
\end{align}
for $x \in [0,1)^d$. If $|y(\nu) - y(0)| >  C_{\rm int} \eps$  we set $\bar{y}_\nu (D^\nu x) = y(0)$. Let $S(\bar{y}_\nu)$ be the set of discontinuity points of $\bar{y}_\nu$ and denote by $\partial_\nu Q^{\nu}$ the two faces of $\partial Q^{\nu}$ which are not parallel to $\nu$. Then $\bar{y}_\nu$ is typically discontinuous on $\partial Q^\nu_\eps(\lambda) \setminus \partial_\nu Q^\nu_\eps(\lambda)$, $\lambda \in \Z^d$. This, however, will not affect our analysis. Essentially, we observe that $S(\bar{y}_\nu) \cap \partial_\nu Q^\nu_\eps(\lambda) \neq \emptyset$ can only occur if there is some cube $Q_\eps(\bar{x})$, $\bar{x} \in \bar{\cal C}'_\eps$ such that $Q_\eps(\bar{x}) \cap Q^\nu_\eps (\lambda) \neq \emptyset$. This is due to \eqref{eq: bar int def} and the definition of ${\cal C}'_\eps$. 

For later we compute the ${\cal H}^{d-1}$ volume of $\partial_\nu Q^\nu$. We first choose $\hat{\nu} \in \R^d$ such that $|\hat{\nu}| = 1$ and $\hat{\nu} \cdot v_{i_j} = 0$ for $j=1, \ldots, d-1$, i.e.\ $\hat{\nu}$ is a unit normal vector to $\partial_\nu Q^\nu$. Then set $\bar{\nu} = |\nu \cdot \hat{\nu}| \hat{\nu}$ and obtain 
\begin{align}\label{eq: face volume}
\frac{1}{2}{\cal H}^{d-1}(\partial_\nu Q^\nu) = \frac{|\det (\bar{\nu}, v_{i_1}, \ldots,v_{i_{d-1}}) |} {|\bar{\nu}|} = \frac{|\det D^\nu|}{|\nu \cdot \hat{\nu}|} =  \frac{\det A}{|\nu \cdot \hat{\nu}|}.
\end{align}
     
We recall here some important properties of the interpolation $\tilde{y}$ on cells $Q_\eps(\bar{x})$, $\bar{x} \in {\cal C}'_\eps$ being proved in \cite{Schmidt:2009} for the case $p=2$. The extension of the results to general $p$ are straightforward.

\begin{lemma}\label{lemma: interpolation}
Let $y: \calL_\eps \to \R^d$ a lattice deformation, $\tilde{y}$ the corresponding linear interpolation. Then for every $1 \leq p < \infty$ there are constants $c,C> 0$ such that for every cell $Q=Q_\eps(\bar{x})$, $\bar{x} \in {\cal C}'_\eps \setminus \bar{\cal C}'_\eps$ we have
\begin{itemize}
	\item[(i)] $c\dist^p(\bar{\nabla} y|_Q, \bar{SO}(d)) \leq \frac{1}{|Q|}\int_Q \dist^p(\nabla \tilde{y}, SO(d)) \leq C\dist^p(\bar{\nabla} y|_Q,                       \bar{SO}(d))$
	\item[(ii)] $c|\bar{\nabla} y|_Q|^p \leq \frac{1}{|Q|}\int_Q |\nabla \tilde{y}|^p \leq C|\bar{\nabla} y|_Q|^p$.
\end{itemize}
\end{lemma}

The interpolation $\tilde{y}$ proves useful to show that in the continuum limit the discrete gradient reduces to a classical gradient (again cf.\ \cite{Schmidt:2009}).

\begin{lemma}\label{lemma: gradient}
Let $U \subset \Omega$, $\eps_k \to 0$ and a sequence $y_k: \calL_\eps(U) \to \R^d$ with $\bar{\nabla} y_k \rightharpoonup f$ in $L^p$, $\tilde{y}_k \rightharpoonup y$ in $W^{1,p}$ for some $f \in L^p(U)$, $y \in W^{1,p}(U)$, $1\leq p < \infty$. Assume that $\bar{\cal C}'_{\eps_{k}} = \emptyset$ for all $k$. Then
$$f = \nabla y \cdot Z.$$
\end{lemma}  

The following lemma shows that passing from $\tilde{y}$ to $\bar{y}_\nu$, $\nu \in {\cal V}$, we do not change the limit.

\begin{lemma}\label{lemma: limit change}
Let $U \subset \Omega$, $\eps_k \to 0$ and $y_k: \calL_{\eps_k} (U) \to \R^d$ be a sequence of lattice deformations with $\# \bar{\cal C}'_{\eps_k}   \leq C \eps_k^{-d + 1}$. Let $\tilde{y}_k$, $\bar{y}_{\nu,k}$ be the corresponding interpolations. Then passing to the limit $\eps_k \to 0$ we obtain 
$\tilde{y}_k - \bar{y}_{\nu,k} \to 0$ in measure for all $\nu \in {\cal V}$.
\end{lemma}
 
\Proof Let $\nu \in {\cal V}$ and consider the sequences $\tilde{y}_k$ and $\bar{y}_{\nu,k}$. By ${\cal D}'_{\eps_k}(U) \subset \calL'_{\eps_k}$ we denote the midpoints $\bar{x}$ of all cells being either a broken cell itself or a neighbor of a broken cell, i.e.\ $\bar{x} \in \bar{\cal C}'_{\eps_k}$ or $\bar{x} + \eps_k \nu \in \bar{\cal C}'_{\eps_k} $ for some $\nu \in {\cal V}$. Note that $\# {\cal D}'_{\eps_k}(U) \leq 3^d \# \bar{\cal C}'_{\eps_k}  \leq C \eps_k^{-d + 1}$. 

We first take a cell $Q_{\eps_k}(\bar{x})$ with $\bar{x} \notin {\cal D}'_{\eps_k}(U)$ into account. By construction we have $\tilde{y}_k(\bar{x} + \eps z_i) = \bar{y}_{\nu,k}(\bar{x} + \eps z_i)$ for a suitable $z_i \in \calZ$. As $|\nabla \tilde{y}_k|, |\nabla \bar{y}^\nu_k| \leq C C_{\rm int}$ on $Q_{\eps_k}(\bar{x})$ we deduce that 
$$\sup_{x \in Q_{\eps_k}(\bar{x})} |\tilde{y}_k(x) - \bar{y}_{\nu,k}(x)| \leq C C_{\rm int} \eps_k$$
for all $\bar{x} \notin {\cal D}'_{\eps_k}(U)$. Noting that 
$$  \sum_{\bar{x} \in \bar{\cal D}'_{\eps_k}(U)} |Q_{\eps_k}(\bar{x})| \leq C\eps^d_k \, \# {\cal D}'_{\eps_k}(U) \leq C \eps_k \to 0$$
as $\eps_k \to 0$, we deduce that $\tilde{y}_k - \bar{y}_{\nu,k} \to 0$ in measure for $\eps_k \to 0$. \eop

\subsection{Sets of finite perimeter and SBV functions}\label{sec: sbv}

As mentioned before the deformation $\tilde{y}$ typically lies in the space of \textit{special functions of bounded variation}. We briefly recall the definition and state fundamental compactness and slicing results. For proofs and further properties we refer to \cite{Ambrosio-Fusco-Pallara:2000}. Let $U \subset \R^d$ be bounded. Recall that the space $SBV(U; \R^{ m})$ consists of functions $y \in L^1(U; \R^{ m})$ whose distributional derivative $Dy$ is a finite Radon measure, which splits into an absolutely continuous part with density $\nabla y$ with respect to Lebesgue measure and a singular part $D^j y$ whose Cantor part vanishes and thus is of the form 
$$ D^j y = [y] \otimes \xi_y {\cal H}^{d-1} \lfloor S(y), $$
where $S(y)$ (the `jump part' of $Dy$) is an ${\cal H}^{d-1}$-rectifiable set in $U$, $\xi_y$ is a normal of $S(y)$ and $[y]= y^+ - y^-$ (the `crack opening') with $y^{\pm}$ being the one-sided limits of $y$ at $S(y)$. 

We state a version of Ambrosio's compactness theorem adapted for our purposes: 
\begin{theorem}\label{th: compact}
Let $(y_\eps)_\eps$ be a sequence in $SBV(U,\R^{ m})$ such that
$$\int_U |\nabla y_\eps(x)|^2 \, dx + {\cal H}^{d-1}(S(y_{\eps})) + \Vert y_\eps\Vert_\infty \le C  $$
for some constant $C$ not depending on $\eps$. Then there is a subsequence (not relabeled) and a function $y \in SBV(U;\R^d)$ such that $y_\eps \to y$ in $L^1(U)$, and 
\begin{align}\label{eq: compctconv}
\begin{split}
& \nabla y_\eps \rightharpoonup \nabla y \text{ in } L^2(U), \\
& D^j y_\eps \to D^j y \text{ as Radon measures.}
\end{split}
\end{align}
\end{theorem}
To state the slicing property we first define the set 
$$U^{\nu,s} = \lbrace t \in \R: s + t\nu \in U\rbrace$$ 
for $\nu \in S^{d-1}$ and $s \in \R^d$. For a function $y: U \to \R^d$ we then define $y^{\nu,s} : U^{\nu,s} \to \R^d$ by
$$y^{\nu,s}(t) = y(s + t\nu).$$

\begin{theorem}\label{th: slic}
Let $y \in SBV(U;\R^{ m})$. For all $\nu \in S^{d-1}$ and ${\cal H}^{d-1}$-a.e.\ $s$ in $\Pi^\nu = \lbrace x: x\cdot \nu = 0\rbrace$ the function $y^{\nu,s}$ belongs to $SBV(U^{\nu,s};\R^d)$ with
$$S(y^{\nu,s}) = \lbrace t \in \R: s + t\nu \in S(y)\rbrace.$$ Moreover, one has
\begin{align*}
\int_{\Pi^\nu} \# S(y^{\nu,s}) \, d{\cal H}^{d-1}(s) = \int_{S(y)} |\xi_y \cdot \nu| \, d{\cal H}^{d-1}.
 \end{align*}
\end{theorem}

An important subset of $SBV$ is given by the indicator functions $\chi_W$, where $W\subset U$ is measurable with ${\cal H}^{d-1}(\partial W)< \infty$. Sets of this form are called \textit{sets of finite perimeter} (cf.\ \cite{Ambrosio-Fusco-Pallara:2000}). As a direct consequence of Theorem \ref{th: compact} we get the following compactness result.
\begin{theorem}\label{th: set per}
Let $(W_\eps)_\eps \subset U$ be a sequence of measurable sets with ${\cal H}^{d-1}(\partial W_\eps) \le C$ for some constant $C$ independent of $\eps$. Then there is a subsequence (not relabeled) and a measurable set $W$ such that $\chi_{W_\eps} \to \chi_W$ in measure for $\eps \to 0$. 
\end{theorem}

\subsection{An estimate on geodesic distances}

We close this preparatory section with a short lemma about the length of Lipschitz curves in sets $W  \subset \R^{d-1}$ of the form \eqref{eq: W form} introduced below: We estimate geodesic distances and the area swept by curves of given length emanating from a common point in terms of the area and surface of $W$. For this purpose, we define $\dist_W(p,q)$ as the infimum of the length of Lipschitz curves in $W$ connecting the points $p,q \in W$ and let ${\cal H}^m$ denote the $m$-dimensional Hausdorff measure.

\begin{lemma}\label{lemma: curves}
There are constants $C,c,c'>0$ (depending on $D$ and $D'$) such that for all $\tilde{W}, W\subset (0,1)^{d-1}$ of the form \eqref{eq: W form} and $\eps$ small enough the following holds:
\begin{itemize}
\item[(i)] $\dist_{\tilde{W}}(p,q) \le C (1 + \eps^{(s-1)(d-3)})$ for all $p,q \in W$
\item[(ii)] For all $p \in W$, $t \in (0,c\eps^{(s-1)(d-2)})$ one has for $\eps$ small enough 
$${\cal H}^{d-1}(\lbrace q \in \tilde{W}: \dist_{\tilde{W}}(p,q) \le t\rbrace) \ge  c' t \eps^{(1-s)(d-2)}.$$ 
\end{itemize}
\end{lemma} 

\Proof
We cover $\tilde{\Omega} := (0,1)^{d-1}$ up to a set of measure zero with the sets $C_{\eps}(\bar{x}) = \bar{x} + (0, \frac{1}{l})^{d-1}$, where $l = \lceil\frac{1}{\hat{\eps}}\rceil$, $\hat{\eps} = \eps^{1-s}$ and $\bar{x} \in I_\eps(\tilde{\Omega}) \subset \frac{1}{l} \Z^{d-1}$. Also set $I_\eps(\tilde{W}) = \{\bar{x} \in I_\eps(\tilde{\Omega}) : \overline{C_\eps(\bar{x})} \subset \tilde{W}\}$. We let $\tilde{V}$ be the connected component of
$$ \bigcup_{\bar{x} \in I_\eps(\tilde{W})} \overline{C_\eps(\bar{x})} \subset \tilde{W}$$ 
with largest Lebesgue measure. We note that for $D'$ sufficiently large $W \subset \tilde{V}$ and thus also $\tilde{V}$ satisfies condition \eqref{eq: W form} possibly passing to a larger $D$. Given $U=\tilde{\Omega}, \tilde{V}$ for two points $p, q \in I_\eps(U)$ we denote the lattice geodesic distance of $p$ and $q$ in $U$, i.e.\ 
the length of the shortest polygonal path $\Gamma_U(p,q) := (x_0 = p, x_1, \ldots , x_n = q)$ with $x_j \in I_\eps(U)$ and $x_{j+1} - x_j = \pm \frac{1}{l}\e_i$ for some $i=1, \ldots, d-1$ connecting $p$ and $q$, by $d_U(p, q)$.

Denote the connected components of $\tilde{\Omega} \setminus \tilde{V}$ by $\tilde{V}_1, \ldots, \tilde{V}_n$ and choose $I_\eps(\tilde{V}_i) \subset I_\eps(\tilde{\Omega})$ such that $\overline{\tilde{V}}_i = \bigcup_{\bar{x} \in I_\eps(\tilde{V}_i)} \overline{C_\eps(\bar{x})}$. It is easy to see that for (i) it suffices to show that $d_{\tilde{V}}(p,q) \le C(1 + \hat{\eps}^{3-d})$ for all $p,q \in I_\eps(\tilde{V})$. Given $p,q \in I_\eps(\tilde{V})$ we first note that $d_{\tilde{\Omega}}(p,q) \le d-1$. Let $\Gamma_{\tilde{\Omega}}(p,q) = (x_0, \ldots, x_m)$ be a (non unique) shortest lattice path connecting $p$ and $q$. If $x_j \in I_\eps(\tilde{W})$ for all $j$ we are finished. Otherwise, for the local nature of the arguments we may assume that $\Gamma_{\tilde{\Omega}}(p,q)$ intersects exactly one $\tilde{V}_i$. Let $x_{j_1}, x_{j_2} \in I_\eps(\tilde{V}_i)$ be the first and the last point in $\tilde{V}_i$, i.e.\ $x_j \notin I_\eps(\tilde{V}_i)$ for $j < j_1$ and $j > j_2$. Then it is elementary to see that 
$$d_{\tilde{V}}(x_{j_1 -1}, x_{j_2 + 1}) \leq C l^{d-3}{\cal H}^{d-2}(\partial \tilde{V}_i) \le C \hat{\eps}^{3-d}{\cal H}^{d-2}(\partial \tilde{V}_i)$$ 
for some $C>0$ not depending on $\tilde{V}$ as the number of cubes at the boundary of $\tilde{V}_i$ can be bounded by $C l^{d-2}{\cal H}^{d-2}(\partial \tilde{V}_i))$. Let $\Gamma_{\tilde{V}}(x_{j_1 -1}, x_{j_2 + 1}) = (y_0, \ldots, y_{\tilde{m}})$ be a shortest path. Then 
$$(x_0, \ldots, x_{j_1-1},y_1,\ldots,y_{\tilde{m}-1},x_{j_2 +1}, \ldots, x_m)$$
is a lattice path in $\tilde{V}$ connecting $p$ and $q$ which shows that $d_{\tilde{V}}(p,q) \le m + \tilde{m}  \le C +  C\hat{\eps}^{3-d}{\cal H}^{d-2}(\partial \tilde{V}_i) \le C +  CD\eps^{(s-1)(d-3)}$. 

To show (ii) we let $p \in W$ and $t \in (0,c\eps^{(s-1)(d-2)})$ for some small $c>0$. Without restriction we may assume $p \in I_\eps(\tilde{V})$. If $\dist_{\tilde{W}}(p,q) \le t$ for all $q \in \tilde{V}$ the assertion is clear. Otherwise, there is some $q \in I_\eps(\tilde{V})$ with $d_{\tilde{V}}(p,q) \ge \dist_{\tilde{V}}(p,q) \ge \dist_{\tilde{W}}(p,q) > \frac{t}{2}$ and a corresponding shortest path $\Gamma_{\tilde{V}}(p,q) = (x_0= p,x_1,\ldots,x_m=q)$ with $x_i \neq x_j$ for $i \neq j$ and $m \ge \bar{m} := \lceil\frac{lt}{2}\rceil$. Now let $U = \bigcup^{\bar{m}}_{j=0} \overline{C_\eps(x_j)}$. Then it is not hard to see that $U \subset \lbrace q \in \tilde{V}: \dist_{\tilde{V}}(p,q) \le t\rbrace$ for $\eps$ small enough and ${\cal H}^{d-1}(U) \ge l^{1-d} \cdot  c' lt \ge  c'\eps^{(1-s)(d-2)} t$, as desired. \eop

\section{Estimates on a mesoscopic cell }\label{sec: mesoscopic cell}

\subsection{Mesoscopic localization}\label{sec: meso-loc}

The goal of this section is the derivation of a lower comparison potential on `large cells' defined on a mesoscopic scale $\eps^s$ with
$$s=\frac{3d-1}{3d}.$$ 
We define the domain $U_\eps = (0,\eps^s\lambda) \times  (0,\eps^{s})^{d-1} $ for $\lambda \in [\lambda_0,2\lambda_0]$, where $\lambda_0 \geq L=L( \sqrt{A^T A},W_{\rm cell},1,\ldots,1)$. We consider the Bravais-lattice defined in Section \ref{sec: model} with a possible translation. For $\rho \in A[0,1)^d$ set $\calL_{\eps,\rho} = \eps\rho + \calL_\eps$.
Let $D, D' > 0$ and suppose that $\tilde{W} \subset (0,1)^{d-1}$ is connected with 
\begin{align}\label{eq: W form}
\begin{split}
\frac{1}{2} < {\cal H}^{d-1}(\tilde{W}) \leq 1, \ \ \ \ {\cal H}^{d-2}(\partial \tilde{W}) < D, 
\end{split}
\end{align}
such that the connected component $W$ of $\lbrace x \in \tilde{W}: \dist(x,\partial \tilde{W}) > D'\eps^{1-s}\rbrace $
with largest Lebesgue measure also satisfies \eqref{eq: W form} for $\eps$ small enough. In what follows, the set $\lbrace 0, \eps^s\lambda\rbrace \times \eps^s W$ will denote the part of the boundary of the mesoscopic cell where we can control the boundary values (see the proof of Lemma \ref{lemma: cr} below, in particular \eqref{eq: difference}ff.). In the proof of Theorem \ref{th: cleavage} we will consider sets $W$ of the above form with the property that ${\cal H}^{d-1}((0,1)^{d-1} \setminus W)$ is small, i.e. ${\cal H}^{d-1}(W) \approx 1$.

We define
\begin{align}\label{eq: boundary cells}
\partial_W (\calL'_{\eps,\rho}(U_\eps)) = \lbrace \bar{x} \in \calL'_{\eps,\rho}:  \ \overline{Q}_\eps(\bar{x}) \cap (\lbrace 0, \eps^s\lambda \rbrace \times \eps^s W)   \neq \emptyset \rbrace.
\end{align} 
Let $y: \calL_{\eps,\rho}(U_\eps) \to \R^d$ be the lattice deformations on $U_\eps$ with corresponding energy
$$E(U_\eps, y) = \eps^{d-1} \sum_{\bar{x} \in { (\calL'_{\eps,\rho} (U_\eps))^\circ}} W_{\rm cell} (\bar{\nabla} y (\bar{x})) + \frac{1}{2}\eps^{d-1} \sum_{\bar{x} \in { \partial_W(\calL'_{\eps,\rho} (U_\eps))}} W_{\rm cell} (\bar{\nabla} y (\bar{x})).$$
The factor $\frac{1}{2}$ takes account of the fact that in the proof of Theorem \ref{th: cleavage} half of the energy of the boundary cells will be assigned to each of the two adjacent mesoscopic cells. Let $\tilde{y}$ denote the interpolation for $y$ defined in Section \ref{sec: interpol}. For $r \in \R$ we will investigate the minimization problem of finding $\inf E(U_\eps, y)$ under certain boundary conditions given as follows. We define the averaged boundary conditions
\begin{align}\label{eq: averaged boundary}\frac{1}{\eps^{s(d-1)}|W|} \int_{\eps^{s}W} \Big(\tilde{y}^1 (\lambda\eps^s,x') - \tilde{y}^1 (0,x')\Big) \, dx' = \lambda\eps^s(1 +  r).
\end{align} 
Here, $x' = (x_2,\ldots,x_d)$, $|W| = {\cal H}^{d-1}(W)$ and $\tilde{y}^1$ denotes the first component of $\tilde{y}$. Moreover, we introduce the condition
\begin{align}\label{eq: broken cell condition}
\bigcup_{\bar{x} \in \bar{\cal B}'_\eps} Q_\eps(\bar{x}) \cap \big( \lbrace0,\eps^s\lambda\rbrace \times \eps^{s} \tilde{W} \big) = \emptyset,
\end{align}
where $\tilde{W}$ is of the form \eqref{eq: W form} and similar as in Section \ref{sec: interpol}, $\bar{\cal B}'_{\eps} \subset \calL'_{\eps,\rho}(U_\eps)$ denotes the cells, where $\text{diam}_1\left\{\bar{\nabla} y (\bar{x}), {\cal Z}\right\} > C^*_{\rm int}$ for a fixed $C^{*}_{\rm int} > 0$ (which may differ from $C_{\rm int}$ to be chosen later).  Note that in contrast to the definition of $\bar{ \cal C}'_\eps$ we consider $\text{diam}_1 \lbrace \cdot , {\cal Z} \rbrace$ only instead of $\text{diam}\lbrace \cdot , {\cal Z} \rbrace$.

We now concern ourselves with the minimization problem
\begin{align*}
M(U_\eps,r) & := \inf\Big\{E(U_\eps, y): y \text{ satisfies } \eqref{eq: averaged boundary} \text{ and } \eqref{eq: broken cell condition}\Big\}.
\end{align*}
Before we state the main theorem of this section we briefly note that $M(U_\eps,r) = 0$ for $-2\leq r \leq 0$ as the averaged boundary conditions may be satisfied by a suitable rotation of the specimen, i.e.\ $y(x)= Rx$ for $R \in SO(d)$.

\begin{theorem}\label{th: meso}
Let $\lambda_0 \geq L$ and let $C_1, C_2 > 0$ be sufficiently large. Let $\delta > 0$ small. Then for all $W \subset (0,1)^{d-1}$ as in \eqref{eq: W form} and $\rho \in A[0,1)^d$ there is a function 
$$ f: \R \times [\eps^s\lambda_0,2\eps^s\lambda_0] \to \R, \qquad
   (r, \lambda) \mapsto f(r, \lambda), $$
which for $r \leq \max\lbrace \eps^{(s-1)(d-3)}, 1\rbrace C_2$ is convex in $r$ and linear in $\lambda$ and for $\eps$ small enough, independently of $\rho$ and $W$, satisfies 
\begin{itemize} 
\item for $r \in \R$, $\lambda \in [\eps^s\lambda_0,2\eps^s\lambda_0]$: $f(r,\lambda) \leq M(U_\eps,r)$, 
\item for $r \in [0,C_1 \sqrt{\eps}]$, $\lambda \in [\eps^s\lambda_0,2\eps^s\lambda_0]$: 
\begin{align}\label{eq: sharp1}
\begin{split}
f(r,\lambda) & = \eps^{s(d-1)}\lambda \omega(|W|) \ \Big(\frac{\alpha_A}{2 \det A} \frac{r^2}{\eps} - \delta\Big) \\ 
&\leq M(U_\eps,r)  \leq \frac{1}{\omega(|W|)} f(r,\lambda) + 4\eps^{sd}\lambda_0\delta, 
\end{split}
\end{align}
\item for $r\geq \max\lbrace \eps^{(s-1)(d-3)}, 1\rbrace C_2 $, $\lambda \in [\eps^s\lambda_0,2\eps^s\lambda_0]$:
\begin{align}\label{eq: sharp2}
\begin{split}
f(r,\lambda)& = \eps^{s(d-1)} \omega(|W|) \ \Big(\frac{\beta_A}{\det A} - \delta\Big)  \\ 
&\leq M(U_\eps,r) \leq  \frac{1}{\omega(|W|)} f(r , \lambda) + 2\eps^{s(d-1)}\delta 
\end{split}
\end{align}
\end{itemize}
for a continuous function $\omega: [0,1] \to \R$ with $\omega(1) = 1$.
\end{theorem} 
The theorem shows that $f(r, \lambda)$ is a lower bound for $M(U_{\eps},r)$ which becomes sharp in the regimes $r \in [0,C_1 \sqrt{\eps}]$ and $ r \in [\max\lbrace \eps^{(s-1)(d-3)}, 1\rbrace C_2,\infty)$ provided that $|W| \approx 1$.

The proof of Theorem \ref{th: meso} is essentially divided into three steps each of which dealing with one particular regime: The elastic regime (Lemma \ref{lemma: el}), the fracture regime (Lemma \ref{lemma: cr}) and the one in between (Lemma \ref{lemma: med}). In addition, in the case $d\ge4$ we need a short additional argument in the intermediate regime  (Lemma \ref{lemma: med2}). It will be convenient to rescale the system in order to obtain a problem on a macroscopic domain not depending on $\eps$. Therefore, we let $\hat{\eps} = \eps^{1-s} = \eps^{1/3d}$, $\hat{U} = (0,\lambda) \times (0,1)^{d-1}$, $\hat{y}: \calL_{\hat{\eps},\rho}(\hat{U}) \to \R^d$ and
\begin{align}\label{eq: rescaled energy}
\begin{split}
\hat{E}(\hat{U}, \hat{y}) 
&= \hat{\eps}^{\,d} \eps^{-1} \sum_{\bar{x} \in (\calL'_{\hat{\eps},\rho}(\hat{U}))^\circ} W_{\rm cell} (\bar{\nabla} \hat{y} (\bar{x})) + \frac{1}{2}\hat{\eps}^{\,d} \eps^{-1} \sum_{\bar{x} \in \partial_W(\calL'_{\hat{\eps},\rho}(\hat{U}))} W_{\rm cell} (\bar{\nabla} \hat{y} (\bar{x})) \\
& = \eps^{-\frac{2}{3}}\sum_{\bar{x} \in (\calL'_{\hat{\eps},\rho}(\hat{U}))^\circ} W_{\rm cell} (\bar{\nabla} \hat{y} (\bar{x})) + \frac{1}{2}\eps^{-\frac{2}{3}}\sum_{\bar{x} \in \partial_W(\calL'_{\hat{\eps},\rho}(\hat{U}))} W_{\rm cell} (\bar{\nabla} \hat{y} (\bar{x})), 
\end{split}
\end{align}
where $\bar{\nabla} \hat{y}$ is defined as in \eqref{eq: discrete grad} replacing $\eps$ by $\hat{\eps}$ and $\partial_W(\calL'_{\hat{\eps},\rho}(\hat{U}))$ as in \eqref{eq: boundary cells} replacing $\eps^s$ by $1$. The averaged boundary conditions now become
\begin{align}\label{eq: constraint}
\dashint_{W} \Big(\tilde{\hat{y}}^1 (\lambda,x') - \tilde{\hat{y}}^1 (0,x')\Big) \, dx' = \lambda(1 + r),
\end{align}
where $\dashint$ denotes the averaged integral and the interpolation $\tilde{\hat{y}}$ is defined as in Section \ref{sec: interpol}. The condition \eqref{eq: broken cell condition} on the boundary cells reads as 
\begin{align}\label{eq: broken cell condition2}
\bigcup_{\bar{x} \in \bar{\cal B}'_{\hat{\eps}}} Q_{\hat{\eps}}(\bar{x}) \cap \big( \lbrace0,\lambda\rbrace \times \tilde{W} \big) 
= \emptyset
\end{align}
and the minimum problem in the rescaled variant becomes
\begin{align}\label{eq: min rescale}
\hat{M}(\hat{U},r) := \inf\left\{\hat{E}(\hat{U}, \hat{y}): \hat{y} \text{ satisfies }   \eqref{eq: constraint}, \eqref{eq: broken cell condition2}\right\}.
\end{align}
It is not hard to see that
\begin{align}\label{eq: scaling}
 M(U_\eps,r)  = \eps^{sd}\hat{M}(\hat{U},r) = \eps^{\frac{3d-1}{3}}\hat{M}(\hat{U},r).
\end{align}

\subsection{Estimates in the elastic regime}\label{sec: meso-elast}

We first determine $\hat{M}(\hat{U},r)$ for $r$ near zero.

\begin{lemma}\label{lemma: el}
Let $C_{\rm el}>0$. For $0\leq r \leq  C_{\rm el} \sqrt{\eps}$ the minimizing problem \eqref{eq: min rescale} satisfies
\begin{align}\label{eq: min est}
|W|\frac{\lambda\alpha_A}{2 \det A} \frac{r^2}{\eps} + o(1) \leq  \hat{M}(\hat{U},r) \leq \frac{\lambda\alpha_A}{2 \det A} \frac{r^2}{\eps} + o(1)
\end{align}
for $\eps \to 0$ with $\alpha_A$ as in \eqref{eq: alpha lax}. Here $o(1)$ is independent of $\rho \in A[0,1)^d$, $W \subset (0,1)^{d-1}$ and $\lambda \in [\lambda_0,2\lambda_0]$ and depends only on $C_{\rm el}$.
\end{lemma}
\Proof In the following we drop the superscript $\hat{\cdot}$ if no confusion arises. We first show that 
\begin{align}\label{eq: first lower bound}
M(U,r) \geq |W| \frac{\lambda \alpha_A}{2 \det A} \frac{r^2}{\eps} + o(1)
\end{align}
for $\eps \to 0$, where $o(1)$ only depends on $C_{\rm el}$. We argue by contradiction. If the claim were false, there would exist a $\delta > 0$, sequences $\eps_k \to 0$, $C_{\rm el} \sqrt{\eps_k} \geq r_k \to 0$, $\lambda_k \in [\lambda_0,2\lambda_0]$, $\rho_k \in A[0,1)^d$, $W_k \subset (0,1)^{d-1}$ satisfying \eqref{eq: W form} as well as a sequence $y_k : \calL_{\hat{\eps}_k,\rho_k}(U_k) \to \R^d$ satisfying \eqref{eq: constraint} with respect to $r_k$, \eqref{eq: broken cell condition2} and $E(U_k,y_k) \leq M(U_k,r_k) + \frac{1}{k}$ such that
\begin{align}\label{eq: false assu}
E(U_k,y_k) \leq |W_k| \frac{\lambda_k \alpha_A}{2 \det A} \frac{r_k^2}{\eps_k} - \delta.
\end{align}
Passing to a subsequence we may assume that $\rho_k \to \rho \in A[0,1)^d$, $\lambda_k \to \lambda \in [\lambda_0,2\lambda_0]$ and $\eps_k^{-\frac{1}{2}} r_k \to r \geq 0$. Moreover, as ${\cal H}^{d-2}( \partial W_k)$ is uniformly bounded in $k$, we may assume that $\chi_{W_k} \to \chi_{W}$ in measure for some $W \subset (0,1)^{d-1}$ with $\frac{1}{2} \leq |W| \leq 1$ by Theorem \ref{th: set per}. As discussed in Section \ref{sec: bc}, there is an obvious choice for an elastic deformation, namely $y^*_{k}(x) = (1 + r_k) x$ for all $x \in \calL_{\hat{\eps}_k,\rho_k}(U_k)$. We note that $\tilde{y}^*_{k}$ satisfies \eqref{eq: constraint} as this interpolation by construction is equal to the linear map $y^*_{k}$. It is elementary to see that
\begin{align}\label{eq: comparison}
\begin{split}
E(U_k, y^*_{k}) =  \frac{(1  + O(\hat{\eps}_k)) \, \lambda_k}{\det A} \frac{1}{\eps_k} \, W_{\rm cell}((1+ r_k) Z) & \leq \frac{C}{\eps_k} \big(Q_{\rm cell}(r_k Z) + o(r_k^2)\big) \\
& \leq \frac{C (C_{\rm el}\sqrt{\eps_k})^{2}}{\eps_k} \leq C.
\end{split}
\end{align} 
As usual we follow the convention of denoting different constants with the same letter. As by \eqref{eq: least energy} and \eqref {eq: rescaled energy} a broken cell contributes an energy of order $\eps_k^{-\frac{2}{3}}$ the comparison with \eqref{eq: comparison} yields 
\begin{align}\label{eq: brokcell cond}
\bar{\cal C}'_{\hat{\eps}_k} = \emptyset
\end{align}
for $k \in \N$ sufficiently large. This together with \eqref{eq: broken cell condition2} shows that $\tilde{y}_k$ is a continuous, piecewise linear interpolation on the set $V_k$, where 
\begin{align}\label{eq: Vk}
V^\circ_k = \bigcup_{\bar{x} \in (\calL'_{\hat{\eps}_k,\rho_k}(U_k))^\circ} Q_{\hat{\eps}_k}(\bar{x}), \ \ \ V_{k} =  V^\circ_k \cup  \bigcup_{\bar{x} \in \partial_{W_k}(\calL'_{\hat{\eps}_k,\rho_k}(U_k))} Q_{\hat{\eps}_k}(\bar{x}). 
\end{align}
Applying Lemma \ref{lemma: interpolation}(i) and Lemma \ref{lemma: energy well} we obtain
\begin{align*}
\int_{V_k } \dist^2(\nabla \tilde{y}_k,SO(d)) & \leq C \int_{ V_k }\dist^2(\bar{\nabla} y_k(\bar{x}), \bar{SO}(d)) \leq C \eps_k \, E(U_k, y_k) \\ 
& \leq C \eps_k\, (M(U_k,r_k) + k^{-1})  \leq C \eps_k \, (E(U_k, y^*_{k}) + k^{-1})  \leq C \eps_k.
\end{align*}
In order to estimate $\int_{V_k} \dist^2(\nabla \tilde{y}_k,SO(d))$ we use the following geometric rigidity result proved in \cite{FrieseckeJamesMueller:02}. If $U \subset \R^d$ is a (connected) Lipschitz domain and $1 < p < \infty$, then there exists a constant $C = C(U,p)$ such that for any $f \in W^{1,p}(U,\R^d)$ there is a rotation $R \in SO(d)$ such that
\begin{align}\label{eq: geometric rigidity}
\left\|\nabla f - R\right\|_{L^p(U)} \leq C \left\|\dist(\nabla f, SO(d))\right\|_{L^p(U)}. 
\end{align}
Therefore, there are a constant $C=C(\lambda_0)$ and rotations $R_k \in SO(d)$ such that
\begin{align}\label{eq: compactness1}
\left\|\nabla \tilde{y}_k - R_k\right\|^2_{L^2(V_k)} \leq C\eps_k.
\end{align}
In fact, $C$ depends only on $\lambda_0$ as all shapes $V_k$ are related to $(0,\lambda_0) \times (0,1)^{d-1} $ through bi-Lipschitzian homeomorphisms with Lipschitz constants of both the homeomorphism itself and its inverse uniformly bounded in $k$, cf.\ \cite{FrieseckeJamesMueller:02}. Up to a not relabeled subsequence we may assume that $R_k \to R$ for some $R \in SO(d)$. By Poincar\'e's inequality we obtain
\begin{align}\label{eq: compactness2}
\left\|\tilde{y}_k - (R_k \,x + c_k)\right\|^2_{L^2(V_k)} \leq C\eps_k
\end{align}
for suitable constants $c_k \in \R^d$. We let $u_k = \frac{1}{\sqrt{\eps_k}} (y_k - (R_k \, \cdot + c_k))$ and obtain by Lemma \ref{lemma: interpolation}(ii)
\begin{align}\label{eq: compactness3}
\left\|\bar{\nabla} u_k\right\|^2_{L^2(V_k)} \leq \frac{C}{\eps_k} \left\|\nabla \tilde{y}_k - R_k\right\|^2_{L^2(V_k)} \leq C.
\end{align}
From \eqref{eq: compactness1} and \eqref{eq: compactness2} we deduce that for a suitable subsequence (not relabeled) $\chi_{V_k} \tilde{u}_k \rightharpoonup \chi_U u$ and $\chi_{V_k} \nabla \tilde{u}_k \rightharpoonup \chi_U \nabla u$ in $L^2$ for some $u \in H^1(U,\R^d)$, where $U = (0,\lambda) \times (0,1)^{d-1}$. Then by \eqref{eq: compactness3} and possibly passing to a further subsequence we find some $f \in L^2(U,\R^{d \times 2^d})$ such that $\chi_{V_k} \bar{\nabla} u_{k} \rightharpoonup \chi_U f$ in $L^2$. By \eqref{eq: brokcell cond} and Lemma \ref{lemma: gradient} we obtain $f = \nabla u \cdot Z$, i.e.\ in particular $\chi_{V_k} \bar{\nabla} u_k \rightharpoonup \chi_U \nabla u \cdot Z$ in $L^2$. 

We now concern ourselves with the averaged boundary condition \eqref{eq: constraint} for the displacement fields $u_k$. We obtain
\begin{align*}
\lambda_k (1 + r_k) & =  \dashint_{W_k} (\tilde{y}^1_k(\lambda_k,x') - \tilde{y}^1_k(0,x')) \, dx' \\
& = (R_k)_{11} \lambda_k + \sqrt{\eps_k }\dashint_{W_k} (\tilde{u}^1_k( \lambda_k,x') - \tilde{u}^1_k(0,x')) \, dx'
\end{align*}
and therefore
\begin{align}\label{eq: limit b}
\dashint_{W_k} (\tilde{u}^1_k(\lambda_k,x') - \tilde{u}^1_k(0,x')) \, dx' = \lambda_k \frac{r_k}{\sqrt{\eps_k}} + \frac{\lambda_k (1- (R_k)_{11})}{\sqrt{\eps_k}}.
\end{align}
We extend $\tilde{u}_k \in H^1(V_k,\R^d)$ to $\hat{u}_k \in H^1(U_k \cup V_k,\R^d)$ such that $\Vert \hat{u}_k\Vert_{H^1(U_k \cup V_k)} \le C \Vert \tilde{u}_k\Vert_{H^1(V_k)}$,  where $C$ may be chosen independently of $\eps$ for the same reasoning as in \eqref{eq: compactness1}. In particular, we note that $\hat{u}^1_k  = \tilde{u}^1_k$ on $\lbrace 0,\lambda_k\rbrace \times W_k$. Defining $\bar{u}_k \in H^1(U,\R^d)$ by $\bar{u}_k (x_1,x'):= \hat{u}_k(\frac{\lambda_k}{\lambda} x_1,x')$ and possibly passing to a further subsequence we may assume $\bar{u}_k \rightharpoonup u$ weakly in $H^1(U)$. Now choosing $r^* \in \R$ such that the trace of $u$ satisfies 
\begin{align}\label{eq: r*}
\dashint_{W} (u^1(\lambda,x') - u^1(0,x')) \, dx = \lambda r^*, 
\end{align}
by the weak continuity of the trace operator, \eqref{eq: limit b} yields 
\begin{align}\label{eq: r-stern}
  r^* = \lim_{k \to \infty} \eps_k^{-\frac{1}{2}}(r_k + 1 - (R_k)_{11}) \in [0, \infty). 
\end{align}
For later use we remark that, since $r_k \ge 0$, the existence of this limit also implies that $R_{11} = 1$. 

We now derive a lower bound for the limiting energy. To this end, note that by Assumption \ref{assu: 1} for $G \in \R^{d \times2^d}$ small we can write $W_{\rm cell} (Z + G) = \frac{1}{2} Q_{\rm cell}(G) +  \eta(G)$ with $\sup\left\{\frac{ \eta(G)}{|G|^2}: |G|\leq \rho\right\} \to 0$ as $\rho \to 0$. Furthermore, let $\chi_k (x):= \chi_{[0,\eps_k^{-1/4}]} (|\bar{\nabla} u_k (x)|)$ and estimate
\begin{align*}
\begin{split}
  E(U_k,y_k) 
  & \ge \hat{\eps}_k^{d} \eps_k^{-1} \sum_{\bar{x} \in (\calL'_{\hat{\eps}_k,\rho_k}(U))^\circ} W_{\rm cell} (R_k \cdot Z + \sqrt{\eps_k} \, \bar{\nabla} u_k (\bar{x})) \\
  & = \hat{\eps}_k^{d} \eps_k^{-1} \sum_{\bar{x} \in (\calL'_{\hat{\eps}_k,\rho_k}(U))^\circ} W_{\rm cell} (Z + \sqrt{\eps_k} \, R_k^{-1}\,\bar{\nabla} u_k (\bar{x})) \\
  & \geq \frac{1}{\eps_k \det A } \int_{V^\circ_k} \chi_k(x) \, W_{\rm cell} (Z + \sqrt{\eps_k} R_k^{-1}\,\bar{\nabla} u_k (x)) \, dx \\
  & \geq \frac{1}{2\det A} \int_{V^\circ_k} \chi_k(x) \Big(Q_{\rm cell}(R_k^{-1} \,\bar{\nabla} u_k(x)) + \eps_k^{-1} \eta( \sqrt{\eps_k} R_k^{-1}\,\bar{\nabla} u_k (x)) \Big)\,dx.
\end{split}
\end{align*}
The second term may be bounded by
$$\chi_k |\bar{\nabla} u_k|^2 \frac{ \eta(\sqrt{\eps_k}\, \bar{\nabla} u_k)}{|\sqrt{\eps_k}\, \bar{\nabla} u_k|^2}. $$
Since $\chi_{V^\circ_k}\chi_k |\bar{\nabla} u_k|^2$ is bounded in $ L^1$ and $\chi_k \frac{\eta(\sqrt{\eps_k}\, \bar{\nabla} u_k)}{|\sqrt{\eps_k}\, \bar{\nabla} u_k|^2} \to 0$ uniformly, we deduce that $\chi_{V^\circ_k} \chi_k |\bar{\nabla} u_k|^2 \frac{\eta(\sqrt{\eps_k}\, \bar{\nabla} u_k)}{|\sqrt{\eps_k}\, \bar{\nabla} u_k|^2}$ converges to zero in $L^1$ as $k \to \infty$. 
Consequently,
$$\liminf_{k \to \infty} E(U_k,y_k) \geq  \liminf_{k \to \infty} \frac{1}{2\det A} \int_{V^\circ_k} Q_{\rm cell}(\chi_k(x) R_k^{-1} \,\bar{\nabla} u_k(x))\,dx. $$
As $R_k \to R$ and $\chi_k \to 1$ boundedly in measure, we obtain $\chi_{V^\circ_k} \chi_k R^{-1}_k  \bar{\nabla} u_k \rightharpoonup \chi_U R^{-1} \nabla u \cdot Z$ weakly in $L^2$ and thus
\begin{align}\label{eq: false assu2}
\begin{split}
\liminf_{k \to \infty} E(U_k,y_k) \geq \frac{1}{2\det A} \int_{U}  Q_{\rm cell}(R^{-1} \,\nabla u(x) \cdot Z)\,dx =: E_{\rm lim}(U,u). 
\end{split}
\end{align}
We now derive a lower bound for $E_{\rm lim}(U,u)$ which will give a contradiction to \eqref{eq: false assu} and thus \eqref{eq: first lower bound} is proved. Since $R_{11} = 1$, we deduce that $R_{1i} = R_{i1} = 0$ for $i=2,\ldots,d$ and therefore $(R^{-1} \, \nabla u(x))_{11} = (\nabla u(x))_{11}$. Applying \eqref{eq: Wred def}, Lemma \ref{lemma: alpha Q} and \eqref{eq: r*} we obtain 
\begin{align}\label{eq: false assu3}
\begin{split}
E_{\rm lim}(U,u) &\geq \frac{1}{2 \det A} \int_U \tilde{Q}((\nabla u(x))_{11}) \, dx \\
& \geq \frac{|W|}{2 \det A} \dashint_{W} \int^\lambda_0 \alpha_A(\partial_1 u^1(x_1,x'))^2 \, dx_1 \, dx' \\
& \geq \frac{|W| \alpha_A}{2 \det A}  \, \lambda \, \Big(\lambda^{-1} \dashint_{W}  \int^\lambda_0 \partial_1 u^1(x_1,x') \, dx_1  \,dx'\Big)^2 \\ & = \frac{|W| \lambda\alpha_A}{2 \det A} (r^*)^2,
\end{split}
\end{align}
where we have used Jensen's inequality. Recalling \eqref{eq: r-stern} we obtain $r^* \geq \lim_{k \to \infty} \frac{r_k}{\sqrt{\eps_k}} = r$ and thus by \eqref{eq: false assu}, \eqref{eq: false assu2} and \eqref{eq: false assu3}
$$ \infty >\frac{|W| \lambda\alpha_A}{2 \det A} r^2 - \delta \geq \liminf_{k \to \infty} E(U_k,y_k) \geq E_{\rm lim}(U,u) \geq \frac{|W|\lambda\alpha_A}{2 \det A} r^2, $$
giving the desired contradiction. 

To see the upper bound in \eqref{eq: min est}, for given $0\leq r \leq C_{\rm el} \, \sqrt{\eps}$, $\rho \in A[0,1)^d$, $\lambda \in [\lambda_0,2\lambda_0]$ we consider the deformation $y: \calL_{\hat{\eps},\rho} \to \R^d$, 
\begin{align}\label{eq: el min}
y(x) = x + \bar{F}(r) \, x
\end{align}
with $\bar{F}(r)$ as in Lemma \ref{lemma: alpha Q}. As $\bar{F}(r)$ depends linearly on $r$ we obtain 
\begin{align*}
E(U,y) &= \frac{\hat{\eps}^{\, d}}{\eps} \sum_{\bar{x} \in (\calL'_{\eps,\rho}(\hat{U}))^\circ} W_{\rm cell} (Z + \bar{F}(r) \cdot Z) + \frac{\hat{\eps}^{\, d}}{2\eps} \sum_{\bar{x} \in \partial_W(\calL'_{\eps,\rho}(\hat{U}))} W_{\rm cell} (Z + \bar{F}(r) \cdot Z) \\
&=  \frac{\lambda (1 + O(\hat{\eps})) }{\det A} \frac{1}{\eps} \,  W_{\rm cell} (Z + \sqrt{\eps}\,\bar{F}(r/\sqrt{\eps}) \cdot Z) \\
&= \frac{\lambda (1 + O(\hat{\eps}))}{2 \det A} \ Q_{\rm cell} (\bar{F}(r/\sqrt{\eps}) \cdot Z) + \frac{\lambda (1 + O(\hat{\eps}))}{2 \det A} \ \frac{ \eta(\sqrt{\eps}\,\bar{F}(r/\sqrt{\eps}))}{\eps},
\end{align*}
where $\eta$ is as before. Note that the second term converges uniformly to zero as $\eps \to 0$ since $\sqrt{\eps}\,\bar{F}(r/\sqrt{\eps}) \leq C r \leq \sqrt{\eps}\, C C_{\rm el}$. Thus, we obtain
$$E(U,y) = \frac{\lambda\alpha_A}{2 \det A} \frac{r^2}{\eps} + o(1).$$ 
 Since $\tilde{y} = y$ satisfies \eqref{eq: constraint} and \eqref{eq: broken cell condition2}, this shows that $M(U,r) \le \frac{\lambda\alpha_A}{2 \det A} \frac{r^2}{\eps} + o(1)$. \eop

\subsection{Estimates in the intermediate regime}\label{sec: meso-intermediate}

We now determine $\hat{M}(\hat{U},r)$ in an intermediate regime. 
\begin{lemma}\label{lemma: med}
Let $C_{\rm med,2} > 0$ arbitrary, $C_{\rm med,1}> 0$ sufficiently large, $1 < p < \frac{4}{3}$. Then there is a constant $C > 0$ such that the minimizing problem \eqref{eq: min rescale} satisfies
\begin{align}\label{eq: min est3}
 \hat{M}(\hat{U},r) \geq C |W|  \lambda \eps^{-\frac{p}{2}} r^p
\end{align}
for $\sqrt{\eps}\, C_{\rm med,1} \leq r \leq C_{\rm med,2}$ as $\eps \to 0$. The constant $C$ is independent of $\rho \in A[0,1)^d$, $W \subset (0,1)^{d-1}$ and $\lambda \in [\lambda_0,2\lambda_0]$. 
\end{lemma}
Note that we only provide a lower bound which might not be sharp.\\

\Proof We follow the previous proof and only indicate the necessary changes. We again drop the superscript $\hat{\cdot}$ if no confusion arises. By Lemma \ref{lemma: energy well} for a suitable constant $c > 0$ and some $1 < p < \frac{4}{3}$ the cell energy $W_{\rm cell}$ may be bounded from below by a function of the form
\begin{align*}
V_\eps(G) = \begin{cases} \eps^{1 - \frac{p}{2}} \,c \, \chi_{\{\dist(G,\bar{SO}(d))\geq \sqrt{\eps}\}} \, \dist^p(G,\bar{SO}(d)) & G \perp V_0, \ |G| \leq C_{\rm int}, \\ 
W_{\rm cell}(G) & \text{else}. \end{cases}
\end{align*}
Then $\eps^{\frac{p}{2}} r^{-p} E(U,y) \geq  {\cal E}(U,y;r)$, where
\begin{align}\label{eq: calE,E}
{\cal E}(U,y;r) := \hat{\eps}^d \eps^{\frac{p}{2} - 1} r^{-p} \Big(\sum_{\bar{x} \in (\calL'_{\hat{\eps},\rho}(U))^\circ} V_\eps(\bar{\nabla} y (\bar{x})) + \frac{1}{2} \sum_{\bar{x} \in \partial_W(\calL'_{\hat{\eps},\rho}(U))} V_\eps(\bar{\nabla} y (\bar{x})) \Big).
\end{align}
We also note that
\begin{align}\label{eq: V}
V_\eps(G) \geq \eps^{1 - \frac{p}{2}} \,c \, \Big(\dist^p(G,\bar{SO}(d)) -   \eps^{\frac{p}{2}}\Big)
\end{align}
for $G \in \R^{d \times 2^d}$ with $G \perp V_0$ and $|G| \leq C_{\rm int}$. We show that for sufficiently small $\eps$
\begin{align}\label{eq: min est4}
{\cal M}(U,r) := \inf\left\{{\cal E}(U,y; r):  y \text{ satisfies } \eqref{eq: constraint}, \eqref{eq: broken cell condition2}\right\} \geq C |W|\lambda
\end{align}
for some $C > 0$ and argue again by contradiction. If \eqref{eq: min est4} were false, there would exist sequences $\eps_k \to 0$, $C_{\rm med,2}\geq r_k \geq C_{\rm med,1}\sqrt{\eps_k}$, $\lambda_k \in [\lambda_0,2\lambda_0]$, $\rho_k \in A[0,1)^d$, $W_k \subset (0,1)^{d-1}$ satisfying \eqref{eq: W form} as well as a sequence $y_k : \calL_{\hat{\eps}_k,\rho_k}(U_k) \to \R^d$ satisfying \eqref{eq: constraint} with respect to $r_k$, \eqref{eq: broken cell condition2} and ${\cal E}(U_k,y_k;r_k) \leq {\cal M}(U_k,r_k) + \frac{1}{k}$ such that
\begin{align}\label{eq: false assu-int}
{\cal E}(U_k,y_k; r_k) \leq \frac{|W_k|\lambda_k}{k}. 
\end{align}
As above we assume that $\rho_k \to \rho$, $\lambda_k \to \lambda$, $\chi_{W_k} \to \chi_{W}$ in measure and $r_k \to r$ up to subsequences. Plugging in the obvious choice for an elastic deformation $y^*_{k}(x) = (1 + r_k)\, x$, $x \in \calL_{\hat{\eps}_k,\rho_k}(U_k)$ we see that
\begin{align}\label{eq: test energy}
{\cal E}(U_k, y^*_{k}; r_k) \leq C\hat{\eps}_k^d \eps_k^{\frac{p}{2} - 1} r_k^{-p} \ \hat{\eps}_k^{-d} \ c \, \eps_k^{1 - \frac{p}{2}} r_k^p = C,
\end{align} 
and thus as in \eqref{eq: brokcell cond} we deduce that $\bar{\cal C}'_{\hat{\eps}_k} = \emptyset$ for all $k$ sufficiently large since by \eqref{eq: least energy} and \eqref{eq: calE,E} a broken cell contributes an energy of order $r_k^{-p} \eps_k^{\frac{p}{2}-\frac{2}{3}} \geq C^{-p}_{\rm med, 2} \eps_k^{\frac{p}{2}-\frac{2}{3}} \gg 1$. Similarly as in the previous proof we obtain by using Lemma \ref{lemma: interpolation}(i) and \eqref{eq: V}
\begin{align*}
\int_{V_k} \dist^p(\nabla \tilde{y}_k,SO(d))  & \leq C \int_{V_k}\dist^p(\bar{\nabla} y_k(\bar{x}), \bar{SO}(d)) \\ 
& \leq  C \int_{V_k} \eps_k^{\frac{p}{2} - 1} V_{\eps_k}(\bar{\nabla} y_k(\bar{x})) +   C \, c\, \eps_k^{\frac{p}{2}} |V_k|\\
& \leq  C r_k^p \, {\cal E}(U_k, y^*_k; r_k) +  \frac{C r_k^p}{k} + C \eps_k^{\frac{p}{2}} |V_k|
\end{align*}
and thus together with \eqref{eq: test energy} and $r_k \geq C_{\rm med,1}\sqrt{\eps_k}$
$$\int_{V_k} \dist^p(\nabla \tilde{y}_k,SO(d)) \leq C r_k^p$$
if $C$ is sufficiently large. By geometric rigidity \eqref{eq: geometric rigidity} and Poincar\'e's inequality there are rotations $R_k \in SO(d)$ and constants $c_k \in \R^d$ such that
\begin{align}\label{eq: rigidity}
\left\|\nabla \tilde{y}_k - R_k\right\|^p_{L^p(V_k)} \leq C\int_{V_k} \dist^p(\nabla \tilde{y}_k,SO(d)) \leq C r^p_k
\end{align}
and
$$\left\|\tilde{y}_k - (R_k  \cdot + c_k)\right\|^p_{L^p(V_k)} \leq C r^p_k.$$
Letting $u_k = \frac{1}{r_k} (y_k - (R_k \, \cdot + c_k))$ we obtain by Lemma \ref{lemma: interpolation}(ii)
\begin{align}\label{eq: uk Absch}
\left\|\bar{\nabla} u_k\right\|^p_{L^p(V_k)} \leq \frac{C}{r_k^p} \left\|\nabla \tilde{y}_k - R_k\right\|^p_{L^p(V_k)} \leq C 
\end{align}
and deduce that for a suitable subsequence (not relabeled) $R_k \to R$ and $\chi_{V_k} \tilde{u}_k \rightharpoonup \chi_U u$ for some $R \in SO(d)$ and $u \in H^1(U,\R^d)$, where $U = (0,\lambda) \times (0,1)^{d-1}$. Then as in the previous proof we derive $\chi_{V_k} \bar{\nabla} u_k \rightharpoonup \chi_U \nabla u \cdot Z$ in $L^p$ possibly after extracting a further subsequence.  

As before, in particular applying the argument in \eqref{eq: limit b} and \eqref{eq: r*}, we obtain that the limit function satisfies the constraint
$$\dashint_{W} (u^1(\lambda,x') - u^1(0,x')) \, dx' = \lambda r^* \in  [0,\infty),$$
where $r^*  = \lim_{k \to \infty} \frac{1}{r_k} ( r_k + 1 - (R_k)_{11}) =  1 +  \lim_{k \to \infty}\frac{1 - (R_k)_{11}}{r_k}\geq 1$. Applying \eqref{eq: V}, Lemma \ref{lemma: interpolation}(i), \eqref{eq: rigidity} and \eqref{eq: uk Absch} we compute
\begin{align*}
  &\liminf_{k \to \infty} {\cal E} (U_k, y_k; r_k) + \limsup_{k \to \infty}  \frac{|V_k| c \eps_k^{\frac{p}{2}}}{r_k^{p} \det A} \\ 
  &\geq \liminf_{k \to \infty} \frac{c}{r_k^p \det A } \int_{V_k} \dist^p(\bar{\nabla} y_k(x), \bar{SO}(d))\,dx  \\
  & \geq \liminf_{k \to \infty} \frac{C}{r_k^p} \int_{V_k} \dist^p(\nabla \tilde{y}_k(x), SO(d))\,dx 
    \geq \bar{C} \int_{U} |\nabla u(x)|^p\,dx
\end{align*}
for some $\bar{C} > 0$. We define 
$${\cal E}_{\rm lim}(U,u;r^*) =\bar{ C} \int_{U} |\nabla u(x)|^p\,dx$$
and then the arguments in \eqref{eq: false assu3}, in particular a slicing argument and Jensen's inequality, yield
$${\cal E}_{\rm lim}(U,u;r^*) \geq \bar{C}|W|\lambda (r^*)^p \geq \bar{C} |W|\lambda$$
since $r^* \ge 1$. Consequently, for $C_{\rm med,1}$ sufficiently large (independent of $\eps_k$) we derive
$$\liminf_{k \to \infty} {\cal E} (U_k, y_k; r_k)\geq \frac{\bar{C}}{2} |W| \lambda.$$
In view of \eqref{eq: false assu-int} this gives the desired contradiction. Thus,  \eqref{eq: min est4} holds and then by \eqref{eq: calE,E} the claim \eqref{eq: min est3} follows. \eop

\subsection{Estimates in the fracture regime}\label{sec: meso-frac}

We now determine $\hat{M}(\hat{U},r)$ for large $r$.

\begin{lemma}\label{lemma: cr}
Let $\lambda_0 \geq L$. For $\max\lbrace \hat{\eps}^{3-d},1\rbrace C_{\rm cr} \leq r$, $C_{\rm cr}$ sufficiently large, the minimizing problem \eqref{eq: min rescale} satisfies
\begin{align}\label{eq: min est2}
\eps^{-s} (|W| - C(1-|W|)) \Big(\frac{\beta_A}{\det A} + o(1)\Big) \leq  \hat{M}(\hat{U},r) \leq \eps^{-s}  \Big(\frac{\beta_A}{\det A} + o(1)\Big)
\end{align}
for $\eps \to 0$.  Here $C>0$ and $o(1)$ are independent of $\rho \in A[0,1)^d$, $W \subset (0,1)^{d-1}$ and $\lambda \in [\lambda_0,2\lambda_0]$.
\end{lemma}

\Proof We again drop the superscript $\hat{\cdot}$ if no confusion arises. We first show that 
\begin{align}\label{eq: first lower bound2}
 \eps^s M(U,r) \geq \frac{ |W|  \beta_A}{\det A} - \hat{C}(1-|W|) + o(1)
\end{align}
for $\eps \to 0$ and some fixed $\hat{C}$ large enough. We again argue by contradiction. If the claim were false, there would exist a $\delta > 0$, sequences $\eps_k \to 0$, $\max\lbrace \hat{\eps}_{k}^{3-d},1\rbrace C_{\rm cr} \leq r_k$, $\lambda_k \in [\lambda_0,2\lambda_0]$, $\rho_k \in A[0,1)^d$, $W_k \subset (0,1)^{d-1}$ satisfying \eqref{eq: W form}  as well as a sequence $y_k : \calL_{\hat{\eps}_k,\rho_k}(U_k) \to \R^d$ satisfying \eqref{eq: constraint} with respect to $r_k$, \eqref{eq: broken cell condition2} and $E(U_k,y_k) \leq M(U_k,r_k) + \frac{1}{k}$ such that
\begin{align}\label{eq: false cr assu}
  \eps^s_k E(U_k,y_k) \leq \frac{|W_k|\beta_A}{\det A} - \hat{C}(1-|W_k|) - 2 \delta.
\end{align}
Up to choosing subsequences we may assume that $\rho_k \to \rho \in A[0,1)^d$, $\lambda_k \to \lambda \in [\lambda_0,2\lambda_0]$ and $W_k \to W \subset (0,1)^{d-1}$.

We again derive a first upper bound of the minimal energy, now by testing with $y^*_k(x) = x \, \chi_{\{x_1 \leq \lambda_k/2\}} + (x + \lambda_k r_k \e_1) \, \chi_{\{x_1 > \lambda_k/2\}}$. It is easy to see that only cells intersecting the set $\lbrace\frac{\lambda_k}{2}\rbrace \times (0,1)^{d-1}$ give an energy contribution. As the quantity of these cells scales like $\hat{\eps}_k^{\, -d+1}$, by \eqref{eq: compatible} we obtain 
 $$E(U_k,y^*_k) \leq C  \hat{\eps}_k^{-d + 1} \hat{\eps}_k^d \eps_k^{-1} c_* = C \eps_k^{-s} c_*,$$
where $c_*$ only depends on $\beta(\nu), \nu  \in {\cal V}$. Thus, $ \eps^s_k E(U_k,y^*_k) \leq C$ for all $k \in \N$ and some $C$ large enough. Then, by \eqref{eq: least energy} it is not hard to see that there is some $\tilde{C}$ such that 
\begin{align}\label{eq: crack length}
\# \bar{\cal C}'_{\hat{\eps}_k} \leq \tilde{C} \hat{\eps}_k^{-d+1},
\end{align}
where $\tilde{C}$ can be chosen independently of $C_{\rm int} \geq 1$. We now choose $C_{\rm int}$ large enough (depending on $\delta$ and possibly larger than the fixed $C^*_{\rm int}$) such that for every partition $\calZ = \calZ_1 \cup \ldots \cup \calZ_n$ and $G \in \R^{d \times 2^d}$ with $\text{diam}\left\{G,{\cal Z}_i\right\} \leq  C_{\rm int}$ and $\min_{i,j} d(G; \calZ_i,\calZ_j) \geq 2^{-2d} C_{\rm int}$ (see Section \ref{sec: interpol}) we obtain (cf.\ Assumption \ref{assu: 1}(iv))
\begin{align}\label{eq: crack energy}
W_{\rm cell} (G) - \sum^n_{i=1} W^{\calZ_i}(G[\calZ_i]) \geq \frac{1}{2}\sum_{1 \leq i, j \leq n \atop i \ne j} \sum_{z_s \in \calZ_i} \sum_{z_t \in \calZ_j} \beta(z_s, z_t) - \frac{\delta}{\tilde{C}}.
\end{align} 
We obtain from Lemma \ref{lemma: energy well} and Lemma \ref{lemma: interpolation}(i),
\begin{align*}
\eps_k^s  E(U_k,y_k) 
  &\ge \frac{1}{2} \hat{\eps}_k^{d-1} \sum_{\bar{x} \in {\cal F}'_{{\hat{\eps}}_k} \setminus \bar{\cal C}'_{{\hat{\eps}}_k}} W_{\rm cell}(\bar{\nabla} y_k(\bar{x})) 
  + \hat{\eps}_k^{d-1} \sum_{\bar{x} \in \bar{\cal C}'_{{\hat{\eps}}_k}} W_{\rm cell}(\bar{\nabla} y_k(\bar{x})) \\
  & \geq C\hat{\eps}_k^{-1} \sum_{\bar{x} \in {\cal F}'_{{\hat{\eps}}_k} \setminus \bar{\cal C}'_{{\hat{\eps}}_k}}\int_{Q_{\hat{\eps}_k}(\bar{x})} \dist^2(\nabla \tilde{y}_k, SO(d)), 
\end{align*}
where ${\cal F}'_{{\hat{\eps}}_k} =  {\cal C}'_{{\hat{\eps}}_k} \cup \partial_{W_k}(\calL'_{\hat{\eps}_k,\rho_k}(U_k))$. Note that $ \dist(F,B_{\sqrt{d}}(0)) \le  \dist(F,SO(d))$ for all $F \in \R^{d \times d}$, where $B_{\sqrt{d}}(0) \subset \R^{d \times d}$ denotes the ball centered at $0$ with radius $\sqrt{d}$. Therefore, with $V_k$ as in \eqref{eq: Vk} and recalling the construction of $\tilde{y}$ in Section \ref{sec: interpol} with uniformly bounded $\nabla \tilde{y}$ we derive 
\begin{align}\label{eq: ela en}
\begin{split}
  \int_{V_k}\dist^2(\nabla \tilde{y}_k, B_{\sqrt{d}}(0)) 
  & \le \sum_{\bar{x} \in {\cal F}'_{{\hat{\eps}}_k} \setminus \bar{\cal C}'_{{\hat{\eps}}_k}}\int_{Q_{\hat{\eps}_k}(\bar{x})} \dist^2(\nabla \tilde{y}_k, SO(d))  + C C^2_{\rm int} \hat{\eps}^d_k \# \bar{\cal C}'_{\hat{\eps}_k} \\ & \le C  \hat{\eps}_k \eps_k^s E(U_k,y_k) + C \hat{\eps}_k 
   \le C \hat{\eps}_k \to 0
\end{split}
\end{align} 
for $\eps_k \to 0.$

In the following we only consider the first component $\tilde{w}_k := \tilde{y}^1_k$ of the deformations. For $\eta > 0$ we enlarge the set $(0,\lambda_k) \times W_k$ and define $W^\eta_k = ((-\eta,\lambda_k + \eta) \times W_k) \cup U_k$. We extend $\tilde{w}_k$ to $W^\eta_k$ by $\tilde{w}_k(x_1,x') = \tilde{w}_k(0,x') + x_1  \e_1$ for $-\eta < x_{ 1} \leq 0$ and $\tilde{w}_k(x_1,x') = \tilde{w}_k(\lambda_k,x') + (x_1 - \lambda_k) \e_1$ for $\lambda_k \leq x_{ 1} \leq \lambda_k + \eta$.  We note that $S(\tilde{w}_k) \cap W^\eta_k \subset U_k$ by \eqref{eq: broken cell condition2}. 
Due to the boundary condition \eqref{eq: constraint} there are points $q^1_k \in \lbrace0\rbrace \times W_k$ and $q^{ 2}_k \in \lbrace\lambda_k\rbrace \times W_k$ such that
\begin{align}\label{eq: difference}
|\tilde{w}_k(q^1_k) - \tilde{w}_k(q^2_k)| \geq \lambda_k (1+ \max\lbrace \hat{\eps}^{3-d}, 1\rbrace C_{\rm cr}). 
\end{align} 
Due to Lemma \ref{lemma: curves}(i) and condition \eqref{eq: broken cell condition2} there is a constant $C=C(D)$ such that 
$$\sup\lbrace|\tilde{w}_k(j,s) - \tilde{w}_k(j,t)|: s,t \in W_k\rbrace \leq C(D) (1 + \hat{\eps}^{3-d}) C^*_{\rm int}$$
for $j=0,\lambda_k$. Choosing $C_{\rm cr}$ sufficiently large this together with \eqref{eq: difference} shows 
\begin{align*}
\inf\lbrace|\tilde{w}_k(p) - \tilde{w}_k(q)|: p,q \in W^\eta_k, p \cdot \e_1=0, q \cdot \e_1= \lambda_k\rbrace \geq \frac{\lambda_k \max\lbrace \hat{\eps}^{3-d},1\rbrace C_{\rm cr} }{2}.
\end{align*}
Let $U^\eta_k = (-\eta, \lambda_k + \eta) \times (0,1)^{d-1}$. For $M = \frac{\lambda_0 C_{\rm cr}}{2}$ we now introduce the truncated function $\tilde{u}_k: U^\eta_k \to \R$ defined by
\begin{align*}
\tilde{u}_k(x) &:= \max\big\{ \ \min\lbrace(\tilde{w}_k(x) - \tilde{w}_k(0,x'))  , M \rbrace,  -M \ \big\} 
\end{align*}
on $W^\eta_k$ and zero elsewhere, where $x'=(x_2,\ldots,x_d)$. Then it is not hard to see that $\tilde{u}_k(0,x') = 0$,  $|\tilde{u}_k(\lambda_k,x')|= M$ for  $x' \in W_k \text{ a.e.}$ and thus
\begin{align}\label{eq: discrete distance}
|\tilde{u}_k(x_1,x')| \le \eta, \  |\tilde{u}_k(\lambda_k - x_1,x')|\ge M - \eta  \ \   \text{  for  } x_1 \in (-\eta,0), x' \in W_k \text{ a.e.\ }
\end{align}
Moreover, $|\nabla \tilde{u}_k| \le |\nabla \tilde{w}_k| + |\nabla_{x'} \ \tilde{w}_k(0,x')| \le |\nabla \tilde{w}_k| + CC^*_{\rm int}$ a.e.\ on $W^{\eta}_k$. 
The lattice deformation corresponding to $\tilde{u}_k$ is denoted by $u_k$, i.e.\ $u_k = \tilde{u}_k|_{\calL_{\hat{\eps}_k,\rho_k}(U_k)}$. Keeping in mind that $W^{\eta}_k$ is open, it is elementary to see that by truncation of the function no further discontinuity points arise, i.e.\ $S(\tilde{u}_k) \cap W^\eta_k \subset S(\tilde{w}_k) \cap U_k$. Moreover, by \eqref{eq: W form} we obtain
\begin{align}\label{eq: rest boundary}
{\cal H}^{d-1}(S(\tilde{u}_k) \setminus W^\eta_k) \le 2(D\eta + 1 - |W_k|).
\end{align}
We now show that we can find a weakly converging subsequence of $(\tilde{u}_k)_k$. To see this, we first note that by  \eqref{eq: crack length} and \eqref{eq: rest boundary} there is some $C>0$ such that ${\cal H}^{d-1}(S(\tilde{u}_k)  ) \leq {\cal H}^{d-1}(S(\tilde{w}_k)) + {\cal H}^{d-1}(S(\tilde{u}_k) \setminus W^\eta_k) \leq C$ for all $k \in \N$. Moreover $\Vert\nabla \tilde{u}_k\Vert_\infty  \le CC^*_{\rm int} + \Vert \nabla \tilde{w}_k\Vert_\infty  \le C (C_{\rm int} + C^*_{\rm int})$ and $\Vert \tilde{u}_k\Vert_\infty \leq M$ for all $k \in \N$ by construction. Now applying Theorem \ref{th: compact} we deduce that there is some $u \in SBV(U^\eta)$ such that up to a subsequence (not relabeled) $\tilde{u}_k \to u$ in the sense of \eqref{eq: compctconv} and a.e., where $U^\eta = (-\eta, \lambda + \eta) \times (0,1)^{d-1}$. By \eqref{eq: discrete distance} for $M$ large enough with respect to $\eta$ the limit function satisfies
\begin{align}\label{eq: limit difference}
&\operatorname{ess\,inf}\lbrace|u(p) - u(q)|: p \in (-\eta,0) \times W, \ q \in (\lambda, \lambda + \eta) \times W \rbrace \geq \frac{M}{2}.
\end{align}
Note that the above compactness theorem implies $\Vert\nabla u\Vert_\infty \le C(C_{\rm int} + C^*_{\rm int})$. We now improve this bound by showing that $\Vert\nabla u\Vert_\infty \le T$ for some $T > 0$ large enough independent of $\delta$ (recall that $C_{\rm int}$ may depend on $\delta$). Let $R = \sqrt{d} + CC^*_{\rm int}$ for some $C>0$ sufficiently large. Then \eqref{eq: ela en} yields
\begin{align*}
\int_{V_k} \dist^2(\nabla \tilde{u}_k, B_{R}(0))  
\le \int_{V_k } \dist^2(\nabla \tilde{y}_k,B_{\sqrt{d}}(0)) 
   \leq C   \hat{\eps}_k \to 0, 
\end{align*}
for $\eps_k \to 0$. Consequently, we get
\begin{align*}
\int_{U} \dist^2(\nabla u, B_{R}(0)) & \le \liminf_{k \to \infty} \int_{V_k} \dist^2(\nabla \tilde{u}_k, B_{R}(0))  = 0,
\end{align*} 
where we used the convexity of $\dist^2(\cdot, B_{R}(0))$. Therefore, 
$|\nabla u| \le R$ a.e.\ in $U$ and by the extension of $\tilde{y}_k$ to $W^\eta_k$ we get $|\nabla u| \le CC^*_{\rm int}$ a.e.\ on $U^\eta \setminus U$. Consequently, choosing $T> 0$ sufficiently large we obtain 
\begin{align}\label{eq: rigid}
|\nabla u| \le T \text{ a.e.\ in } U^\eta.
\end{align}
We now concern ourselves with the energy contribution of the broken cells $\bar{\cal C}'_{\hat{\eps}_k}$. For all $\nu \in {\cal V}$ we let $\bar{y}_{\nu,k}: U_k \to \R^d$, $\bar{w}_{\nu,k}:  U_k \to \R$ and $\bar{u}_{\nu,k}: U^\eta_k \to \R$ be the interpolations introduced in Section \ref{sec: interpol}. Recalling the construction of the interpolations we obtain by \eqref{eq: beta nu} and \eqref{eq: crack energy}  
\begin{align*}
  \eps_k^s  E(U_k,y_k) 
  &  \ge \hat{\eps}_k^{d-1} \sum_{\bar{x} \in \bar{\cal C}'_{{ \hat{\eps}}_k}} W_{\rm cell}(\bar{\nabla} y_k(\bar{x})) \\
  & \geq  \sum_{\nu \in {\cal V}}  \int_{S(\bar{y}_{\nu,k}) \cap \Gamma(\nu)} \frac{\hat{\eps}_k^{d-1}}{\frac{1}{2}{\cal H}^{d-1}(\partial_\nu Q^\nu_{\hat{\eps}_k})}\beta(\nu) \,d{\cal H}^{d-1} 
 - \hat{\eps}_k^{d-1} \# \bar{\cal C}'_{\hat{\eps}_k} \frac{\delta}{\tilde{C}},
\end{align*}
where $\Gamma(\nu) = \bigcup_{\lambda \in \Z^d} \partial_\nu Q^\nu_{\hat{\eps}_k}(\lambda)$. Then by \eqref{eq: face volume} and \eqref{eq: crack length} we get
$$\eps_k^s  E(U_k,y_k) \ge \sum_{\nu \in {\cal V}} \int_{S(\bar{y}_{\nu,k})} \frac{\beta(\nu)}{\det A} \, |\nu \cdot \xi_{\bar{y}_{\nu,k}}|\,d{\cal H}^{d-1}  - \delta. $$
Applying \eqref{eq: rest boundary} it is not hard to see that there is some $\Gamma_{k, \nu}$ with ${\cal H}^{d-1}(\Gamma_{k,\nu}) \le C(D\eta+ 1- |W_k|)$ such that $S(\bar{u}_{\nu,k}) \subset S(\bar{w}_{\nu,k}) \cup \Gamma_{k,\nu} \subset S(\bar{y}_{\nu,k}) \cup \Gamma_{k,\nu}$. Furthermore, the normals $\xi_{\bar{y}_{\nu,k}}$ and $\xi_{\bar{u}_{\nu,k}}$ coincide on $S(\bar{u}_{\nu,k}) \cap S(\bar{y}_{\nu,k})$. Therefore, we derive
\begin{align}\label{eq: false cr assu2}
\begin{split}
  \delta +  C(D\eta & + 1 - |W_k|) +  \eps_k^s  E(U_k,y_k) \\ & \ge \sum_{\nu \in {\cal V}} \int_{S(\bar{u}_{\nu,k})}  \frac{\beta(\nu)}{\det A} \, |\nu \cdot \xi_{\bar{u}_{\nu,k}}|\,d{\cal H}^{d-1} =: E_S (U_k,u_k). 
  \end{split}
\end{align}
With the notation introduced in Section \ref{sec: sbv} we get using Theorem \ref{th: slic}
\begin{align*}
E_S (U_k,u_k) \geq \frac{1}{\det A} \sum_{\nu \in {\cal V}}  \int_{\Pi^\nu} \# S(\bar{u}^{\nu,s}_{\nu,k}) \beta(\nu) \, d{\cal H}^{d-1}(s).
\end{align*}
Then by the equiboundedness of $E_S(U_k,u_k)$ and Fatou's lemma we deduce that $\liminf_{k \to \infty} \# S(\bar{u}^{\nu,s}_{\nu,k}) < +\infty$ for a.e.\ $s \in \Pi^\nu$ and all $\nu \in {\cal V}$. As $\bar{u}_{\nu,k}$ and $\nabla \bar{u}_{\nu,k}$ are uniformly bounded, by Theorem \ref{th: compact} and Lemma \ref{lemma: limit change} $\bar{u}^{\nu,s}_{\nu,k}$ converges (up to a subsequence) to $u^{\nu,s}$ in the sense of \eqref{eq: compctconv} for a.e.\ $s \in \Pi^\nu$. In particular, we get 
$$\liminf_{k \to \infty} \# S(u^{\nu,s}_{\nu,k}) \geq \# S(u^{\nu,s}_\nu).$$
Applying Fatou's lemma and the slicing theorem once more we then derive
\begin{align}\label{eq: false cr assu3}
\begin{split}
\liminf_{k \to \infty} E_S (U_k,u_k) &\geq \frac{1}{\det A} \sum_{\nu \in {\cal V}}  \int_{\Pi^\nu} \# S(u^{\nu,s}) \beta(\nu) \, d{\cal H}^{d-1}(s)\\
& = \frac{1}{\det A}  \int_{S(u)} \, \sum_{\nu \in {\cal V}}\beta(\nu) |\nu \cdot \xi_u| \, d{\cal H}^{d-1}  =: E_{S, \rm lim} (U,u). 
\end{split}
\end{align}
By \eqref{eq: beta def} and slicing in $\e_1$-direction we get
\begin{align*}
E_{S, \rm lim} (U,u) 
&\ge \frac{1}{\det A} \int_{S(u)} \beta_A |\e_1 \cdot \xi_u| \, d{\cal H}^{d-1} \\ 
&= \frac{1}{\det A}  \int_{(0,1)^{d-1}} \beta_A  \, \# S(u^{\e_1,s}) \, d{\cal H}^{d-1}(s).
\end{align*}
We now choose $M = \frac{\lambda_0 C_{\rm cr}}{2}$ sufficiently large (independently of $\delta$) such that $M \ge 4T\lambda$. Then due to \eqref{eq: limit difference} and \eqref{eq: rigid} it is not hard to see that $\# S(u^{\e_1,s})\geq 1$ for a.e.\ $s\in W$ and therefore $E_{S, \rm lim} (U,u) \geq \frac{|W| \beta_A}{\det A}$.
Letting $\eta \to 0$ and choosing $\hat{C}$ sufficiently large we now conclude by \eqref{eq: false cr assu},  \eqref{eq: false cr assu2} and  \eqref{eq: false cr assu3}:
\begin{align*}
\infty > \frac{|W|\beta_A}{\det A} - 2 \delta & \geq \liminf_{k \to \infty}\eps_k^s E(U_k, y_k) +  \hat{C}(1- |W|)  \geq \liminf_{k \to \infty} E_S(U_k,y_k) - \delta \\
& \geq E_{S, \rm lim}(U,u) - \delta \geq \frac{|W| \beta_A}{\det A} - \delta.
\end{align*}
This gives the desired contradiction. 

To see the upper bound in \eqref{eq: min est2} we choose $\xi \in S^{d-1}$ such that \eqref{eq: beta def} is minimized and define the hyperplane $\Pi = \lbrace x \in \R^d: x \cdot \xi = c\rbrace$ for a suitable $c \in \R$ such that $\Pi \cap U \subset \lbrace \delta \leq x_1 \leq \lambda - \delta \rbrace$ for some $\delta > 0$. We set
\begin{align}\label{eq: cr min}
y(x) = \begin{cases}  x, &  x\cdot \xi \leq c,\\
                      x + r \lambda \e_1,   & x\cdot \xi  >  c.    \end{cases}
\end{align}
The energy corresponding to the deformation $y$ is given by the bonds intersecting $\Pi$. These springs, associated to the lattice directions $\nu \in {\cal V}$, are elongated by a factor scaling with $r/\hat{\eps}$ and yield a contribution $\beta(\nu)$ in the limit $\eps \to 0$ by \eqref{eq: beta nu}. As the projection in $\nu$-direction onto the hyperplane $\lbrace x \cdot \xi = c\rbrace$ of the face $\partial_\nu Q^\nu$ has ${\cal H}^{d-1}$-volume 
$$\frac{1}{2} {\cal H}^{d-1}(\partial_\nu Q^\nu) \, \Big|\frac{\nu}{|\nu|} \cdot \hat{\nu} \Big| \frac{1}{\big|\frac{\nu}{|\nu|} \cdot \xi \big|} 
= \frac{ \hat{\eps}^{d-1} \, \det A}{|\nu \cdot \hat{\nu}|} \frac{|\nu \cdot \hat{\nu} |}{|\nu \cdot \xi |} 
= \frac{\hat{\eps}^{d-1}\, \det A }{|\nu \cdot \xi |}$$ 
(see \eqref{eq: face volume}) it is not hard so see that 
\begin{align}\label{eq: counting}
\frac{|\nu \cdot \xi|}{\hat{\eps}^{d-1} \det A |e_1 \cdot \xi| }  + O\Big(\frac{1}{\hat{\eps}^{d-2}}\Big)
\end{align}
springs in $\nu$-direction are broken. This yields the energy
$$ \eps^{-s}\Big(\frac{\beta_A}{\det A} + o(1)\Big) + O(\hat{\eps}^{2} \eps^{-1}),$$
for $\eps \to 0$ as desired.  \eop 

\subsection{Estimates in a second intermediate regime}\label{sec: meso-intermediate-d4}

We provide an additional lemma needed in the case $d \ge 4$, which may be safely skipped by a reader who is more interested in the physical application of our result. 

\begin{lemma}\label{lemma: med2}
Let $d \ge4$, $C^*_{\rm med,2} > 0$ arbitrary and $C^*_{\rm med,1}> 0$ sufficiently large. Then there is a constant $C > 0$ such that the minimization problem \eqref{eq: min rescale} satisfies
\begin{align*}
 \hat{M}(\hat{U},r) \geq C \eps^{d-2 + s(1-d)} r
\end{align*}
for $C^*_{\rm med,1} \leq r \leq C^*_{\rm med,2} \hat{\eps}^{3-d}$ as $\eps \to 0$. The constant $C$ is independent of $\rho \in A[0,1)^d$, $\tilde{W} \subset (0,1)^{d-1}$ and $\lambda \in [\lambda_0,2\lambda_0]$.
\end{lemma}

\Proof 
The superscript $\hat{\cdot}$ is again dropped where no confusion arises. Let $C^*_{\rm med,2} > 0$. Let $\rho, \tilde{W}, \lambda$ and $r$ with $C^*_{\rm med,1} \leq r \leq C^*_{\rm med,2} \hat{\eps}^{3-d}$ be given and consider a deformation $y : \calL_{\hat{\eps},\rho}(U) \to \R^d$ satisfying \eqref{eq: constraint} with respect to $r$, \eqref{eq: broken cell condition2} and $E(U,y) \le 2 M(U,r)$. Due to \eqref{eq: constraint}  there is a $q \in W$ such that $|\tilde{y}^{1}(\lambda,q) - \tilde{y}^{1}(0,q)| \ge \lambda(1 + r)$. Applying  Lemma \ref{lemma: curves}(ii) for $t = \frac{\lambda(1 + r)}{4 C C^*_{\rm int}}$ and \eqref{eq: broken cell condition2} we find a set $V\subset \tilde{W}$ with ${\cal H}^{d-1}(V) \ge c' \lambda_0 r \eps^{(1-s)(d-2)}$ such that  $$|\tilde{y}(\lambda,q) - \tilde{y}(0,q)| \ge \frac{\lambda(1 + r)}{2} \ge  \frac{\lambda_0 C^*_{\rm med,1}}{2}$$ 
for all $q \in V$. Fix $C_{\rm int}$ as defined in Section \ref{sec: interpol} and recall $\Vert \nabla \tilde{y}\Vert_\infty \le CC_{\rm int}$. Choose $C^*_{\rm med,1}$ large enough such that $((0,\lambda) \times \lbrace q\rbrace) \cap \bigcup_{\bar{x} \in \bar{\cal C}'_{\hat{\eps}}} \overline{Q_{\hat{\eps}}(\bar{x})} \neq \emptyset$ for all $q \in V$. As the orthogonal projection of a cell onto $\lbrace0\rbrace \times \R^{d-1}$ has ${\cal H}^{d-1}$-measure smaller than $C \hat{\eps}^{d-1}$  we deduce
 $$\# \bar{\cal C}'_{\hat{\eps}} \ge C\eps^{s-1}r$$
 for some $C>0$. As every broken cell provides at least the energy $C\eps^{-\frac{2}{3}} = C\eps^{d(1-s)-1}$ by \eqref{eq: least energy} we derive
 $$M(U,r) \ge C\eps^{d(1-s)-1}\ \# \bar{\cal C}'_{\hat{\eps}} \ge C \eps^{d-2 + s(1-d)}r .$$
\eop

\subsection{Proof of Theorem \ref{th: meso}}\label{sec: pf41}

Summarizing our previous estimates we are now in a position to prove Theorem \ref{th: meso}. 
\smallskip

\noindent {\em Proof of Theorem \ref{th: meso}.}
We begin to construct such a function $\hat{f}: \R \times [\lambda_0,2\lambda_0] \to \R$ for the rescaled problem $\hat{M}(\hat{U},r)$. Let $\omega(|W|) = |W|(1 - C(1 - |W|))$ with the constant $C$ of Lemma \ref{lemma: cr}. For $\delta > 0$ small we set $\hat{f}(r,\lambda) = -\delta  \lambda \omega(|W|) $ for $r\leq 0$ and for $C_1 > 0$ sufficiently large such that $C_1 \ge C_{\rm med, 1}$ we define 
$$\hat{f}(r,\lambda) = \omega(|W|) \lambda \Big(\frac{\alpha_A}{2 \det A} \frac{r^2}{\eps} - \delta\Big)$$
for $r \in [0, C_1 \sqrt{\eps}]$, $\lambda \in [\lambda_0,2\lambda_0]$. Choose the affine function $g: \R \to \R$ satisfying $\lambda \omega(|W|) g(C_1 \sqrt{\eps}) = \hat{f}(C_1 \sqrt{\eps},\lambda)$ and $\lambda \omega(|W|) g' = \partial_r\hat{f}(C_1 \sqrt{\eps},\lambda)$. For $t > C^p_1 C$ sufficiently large we let $h(r) = C \eps^{-\frac{p}{2}}  r^p - t$ for $r\geq C_{\rm med,1} \sqrt{\eps}$ with $C_{\rm med,1}$ as in Lemma \ref{lemma: med}, so that there is a (unique) intersection point of the graphs of $g$ and $h$, $(\bar{r}_t,g(\bar{r}_t)) = (\bar{r}_t,h(\bar{r}_t))$, for which $h'(\bar{r}_t) \geq g'$. Note that $\bar{r}_t \sim \sqrt{\eps}$. Then we set $\hat{f}(r,\lambda) = \lambda \omega(|W|) g(r)$ for $r \in [C_1 \sqrt{\eps},\bar{r}_t]$ and $\hat{f}(r,\lambda) = \lambda\omega(|W|) h(r)$ for $r \in [\bar{r}_t, \max\lbrace \hat{\eps}^{3-d},1\rbrace C_2]$ for $C_2 >0$ so large that $C_2 \ge \max \lbrace C_{\rm cr}, C^*_{\rm med,1}\rbrace$. Finally, we let 
$$\hat{f}(r,\lambda) = \eps^{-s} \omega(|W|) \Big(\frac{\beta_A}{\det A} - \delta\Big)$$
for $r \geq \max\lbrace \hat{\eps}^{3-d},1\rbrace C_2$. In the case $d \ge 4$ we observe that for $C_2\le r \le C_2\hat{\eps}^{3-d}$ we have 
$$ \omega(|W|) \lambda \eps^{-\frac{p}{2}} r^p \le C \eps^{d-2 + s(1-d)} r \le C\eps^{-s} \hat{\eps} \le C\eps^{-s}$$ 
for some $p>1$ small enough. Now by construction it is not hard to see that $\hat{f}$ is convex for $r \leq \max\lbrace \hat{\eps}^{3-d},1\rbrace C_2$. Moreover, we obtain $\hat{f} \leq \hat{M}(\hat{U},\cdot)$ for $\eps$ small enough independently of $\rho \in A[0,1)^d$ and $W$. Indeed, for $r \in [0, \bar{r}_t]$ this follows from Lemma \ref{lemma: el}, for $r \in [\bar{r}_t,C_2]$ from Lemma \ref{lemma: med}, for $r \in [\hat{\eps}^{3-d} C_2,\infty)$ from Lemma  \ref{lemma: cr} and in the case $d \ge 4$ we use Lemma \ref{lemma: med2} for the additional regime $[C_2, \hat{\eps}^{3-d} C_2]$. In particular, observe that $\bar{r}_t \ge C_{\rm med,1}$ and $C_2 \ge \max \lbrace C_{\rm cr}, C^*_{\rm med,1}\rbrace$ such that the lemmas can be applied. Furthermore, we derive 
$$\hat{f}(r,\lambda)  \leq \hat{M}(\hat{U},r) \leq \frac{1}{\omega(|W|)}\hat{f}(r,\lambda) + 4\lambda_0\delta $$
for $r \in [0,C_1 \sqrt{\eps}]$ and
$$\hat{f}(r,\lambda)  \leq \hat{M}(\hat{U},r) \leq \frac{1}{\omega(|W|)}\hat{f}(r,\lambda) + 2 \eps^{-s}\delta$$
for $r \geq \max\lbrace \hat{\eps}^{3-d},1\rbrace C_2$. To finish the proof it suffices to recall $M(U_\eps, r) = \eps^{sd} \hat{M}(\hat{U},r)$ by \eqref{eq: scaling} and to set $f(r,\lambda)= \eps^{sd} \hat{f}(r, \eps^{-s}\lambda)$ for all $r\in \R$ and $\lambda \in \eps^s[\lambda_0,2\lambda_0]$. \eop

\section{Proof of the main theorem}\label{sec: proof}

We are now in a position to prove the main theorem.
\smallskip

\noindent {\em Proof of Theorem \ref{th: cleavage}.}
Let $y \in {\cal A}(a_\eps)$. We partition $(0,l_2) \times \ldots \times (0,l_d)$ up to a set of size $O(\eps^s)$ with sets $V_i$, $i \in I$, which are translates of the cube $\eps^s(0,1)^{d-1}$. Furthermore, we set $V^d_i = (0,l_1) \times V_i$ for all $i \in I$. For $C^*_{\rm int} > 0$ we denote the set of broken cells by $\bar{\cal B}'_\eps$ as defined at the beginning of Section \ref{sec: mesoscopic cell}. We let
$$\bar{I} := \Big\lbrace i \in I:  \# \lbrace\bar{x} \in \bar{\cal B}'_\eps: Q_\eps(\bar{x}) \subset V^d_i\rbrace >  \frac{2\beta_A}{\eps^{(1-s)(d-1)} \bar{C}\det A} \Big\rbrace$$
with $\bar{C} =  \bar{C}(C^*_{\rm int})$ as in \eqref{eq: least energy}. Then for $i \in \bar{I}$ we estimate
\begin{align}\label{eq: broken energy}
  \eps^{d-1} \sum_{\bar{x} \in (\calL'_\eps(V^d_i))^\circ} W_{\rm cell}(\bar{\nabla} y(\bar{x})) 
  \geq \eps^{d-1}  \frac{2\beta_A}{\eps^{(1-s)(d-1)} \bar{C}\det A} \,  \bar{C} 
  = \frac{\eps^{s(d-1)} 2 \beta_A}{\det A}. 
\end{align}
Now consider some $V_i$ for $i \in I \setminus \bar{I}$. For $\lambda_0 \geq L( \sqrt{A^TA},W_{\rm cell},1,\ldots,1)$ we partition $V^d_i$ into sets of the form $U_1 = (u_0,u_1) \times V_i, \ldots, U_n = (u_{n-1},u_n) \times V_i$, where $u_0 = l_A\eps$, $u_n = l_1 - l_A\eps$ and $u_j - u_{j-1} \in \eps^s[\lambda_0,2\lambda_0]$ for all $j=1,\ldots,n$. This can and will be done so that the number of broken cells $Q_\eps(\bar{x}) \subset V_i^d$, $\bar{x} \in \bar{\cal B}'_\eps$, intersecting $\lbrace u_j \rbrace \times V_i$ is minimal on an interval of length $\frac{\lambda_0 \eps^s}{8}$ with boundary point $u_j$, i.e., 
\begin{align}\label{eq: vj cond}
N(u_j) &:= \# {\cal T}(u_j)  = \min_{\bar{u} \in J(u_j)} \#  {\cal T}(\bar{u}),
\end{align}
for all $j=1,\ldots,n-1$, where 
$$  {\cal T}(s) = \lbrace \bar{x} \in \bar{\cal B}'_\eps : 
   Q_\eps(\bar{x}) \subset V_i^d \text{ and } 
   Q_\eps(\bar{x}) \cap (\lbrace s \rbrace \times V_i) \neq \emptyset\rbrace $$ 
and either $J(u_j) = [u_j - \frac{\lambda_0 \eps^s}{8}, u_j]$ or $J(u_j) = [u_j, u_j +  \frac{\lambda_0 \eps^s}{8}]$. Such values of $u_j$ can be constructed by first considering the equidistant points $v_j = l_A \eps + j \mu \eps^s$ with $l_A \eps + n \mu \eps^s = l_1 - l_A \eps$ and $|\mu - \frac{3 \lambda_0}{2}| \le \frac{\lambda_0}{4}$ and then choosing $u_j = \operatorname{arg\,min} \{ \# {\cal T}(\bar{u}) : v_j - \frac{\lambda_0 \eps^s}{8} \le \bar{u} \le v_j + \frac{\lambda_0 \eps^s}{8} \}$ for $j = 1, \ldots, n-1$. 

Moreover, we have $N(u_0) = N(u_n) = 0$ due to the boundary conditions \eqref{eq: bc1}. We now show that 
\begin{align}\label{eq: vj cond2}
  \sum^n_{j=0} N(u_j)
  \leq \frac{C \eps^{(s-1)(d-2)}}{\lambda_0}.
\end{align}
We cover $J(u_j) \times V_i$ with translates of $(0,\eps l_A) \times (0,\eps^s)^{d-1}$, where $l_A$ is as defined in \eqref{eq: lA}. As every cell is contained in at most two of these translates we derive
\begin{align*}
  \# \lbrace\bar{x} \in \bar{\cal B}'_\eps : 
   Q_\eps(\bar{x}) \subset V_i^d \text{ and } Q_\eps(\bar{x}) \cap (J(u_j) \times V_i) \neq \emptyset\rbrace 
  \geq \Big\lfloor\frac{ \lambda_0 \eps^s}{16 l_A \eps}\Big\rfloor N(u_j) 
\end{align*}
for $j=1,\ldots,n-1$ due to the construction \eqref{eq: vj cond}. Summing over $j$, we find 
$$ \sum^n_{j=0} N(u_j) \le \frac{C \eps}{\eps ^s \lambda_0} \# \lbrace\bar{x} \in \bar{\cal B}'_\eps: Q_\eps(\bar{x}) \subset V_i^d \rbrace 
   \le \frac{C \eps^{1 - s} \eps^{(1 - s)(1 - d)}}{\lambda_0} 
   = \frac{C \eps^{(s-1)(d - 2)}}{\lambda_0} $$
since $i \in I \setminus \bar{I}$. Note that the estimate \eqref{eq: vj cond2} relies only on the fact that $i \in I \setminus \bar{I}$ but is independent of the particular set $V_i$, the deformation $y$ and $\eps$.  

Let $T_i = \bigcup^{n-1}_{j=1} \bigcup_{\bar{x} \in {\cal T}(u_j)} \overline{Q_{\eps}(\bar{x})}$ and $S_i = \pi_1 T_i$, where $\pi_1 T_i \subset \R^{d-1}$ denotes the set which arises from $T_i$ by orthogonal projection onto $\lbrace0\rbrace \times V_i$ and cancellation of the first component. Using \eqref{eq: vj cond2} we find 
$${\cal H}^{d-2}(\partial S_i) \le \sum^{n-1}_{j=1} N(u_j) {\cal H}^{d-2}(\partial \, \pi_1 Q_{\eps}) \le C\eps^{d-2}\sum^{n-1}_{j=1} N(u_j) \le C \lambda_0^{-1} \eps^{s(d-2)}.$$
We choose $\lambda_0$ so large that ${\cal H}^{d-2}(\partial S_i) \le \delta\eps^{s(d-2)}$. Let $V_{i, \eps} = \{ x \in V_i : \dist(x, \partial V_i) \ge C \eps\}$ with $C$ so big that $\pi_1 Q_{\eps}(\bar{x}) \cap V_{i, \eps} = \emptyset$ whenever $Q_{\eps}(\bar{x}) \not\subset V_i^d$. By the isoperimetric inequality we deduce that there is a unique connected component $\tilde{W}_i$ of $V_{i,\eps} \setminus S_i$ satisfying $|\tilde{W}_i|:= {\cal H}^{d-1}(\tilde{W}_i) \geq (1 - C \eps^{1 - s} - C\delta^{\frac{d-1}{d-2}}) \eps^{s(d-1)} $, where $C$ is a constant only depending on the dimension. Moreover, we have ${\cal H}^{d-2} (\partial \tilde{W}_i) \le C\eps^{s(d-2)}$ and so we see that for $\delta$ small enough $\tilde{W}_i$ is of the form \eqref{eq: W form} (after rescaling by $\eps^{-s}$). Furthermore, by a similar argument (e.g.\ by enlarging the cubes which form $T_i$) we find that 
$$ {\cal H}^{d-2} \big( \lbrace x \in  \tilde{W}_i: \dist(x,\partial \tilde{W}_i) = D'\eps \text{ and } 
   \dist(x,\partial V_{i, \eps}) \ne D' \eps\rbrace \big) 
   \le C\delta\eps^{s(d-2)}.$$
Consequently, we define $W_i$ corresponding to $\tilde{W}_i$ as described in \eqref{eq: W form} (replacing $\eps^{1-s}$ by $\eps$ due to the different scaling) and obtain $|W_i|:= {\cal H}^{d-1}(W_i) \geq (1 - C \eps^{1 - s} - C\delta^{\frac{d-1}{d-2}}) \eps^{s(d-1)}$ and ${\cal H}^{d-2} (\partial W_i) \le C\eps^{s(d-2)}$ for some possibly larger constant $C$. Clearly, $W_i$ is of the form \eqref{eq: W form}. The sets $U_j$ defined above correspond to $U_\eps$ considered in Section \ref{sec: mesoscopic cell} up to a translation. In particular, the sets $\tilde{W}_i$ satisfy condition \eqref{eq: broken cell condition} due to the construction of $S_i$. 

We define
$$ r_j :=  - 1 + \frac{1}{u_j - u_{j-1}}\dashint_{W_i} \Big(\tilde{y}^1 (u_j,x') - \tilde{y}^1 (u_{j-1},x')\Big) \, dx' $$ 
for $j=1, \ldots,n$. Note that this definition is meaningful as $\tilde{y}^1$ is defined on all of $(0, l_1) \times W_i$ (see Section \ref{sec: interpol}). As $y \in {\cal A}(a_\eps)$ it is not hard to see that 
\begin{align*}
\sum^n_{j=1}  (u_j - u_{j-1}) \, r_j & = -(l_1  -2l_A\eps) + \dashint_{W_i} \Big(\tilde{y}^1 (l_1 - l_A\eps,x') - \tilde{y}^1 (l_A\eps,x')\Big) \, dx'\\ &  = - (l_1  -2l_A\eps) + l_1 - 2l_A \eps + (l_1 -2l_A\eps)a_\eps = (l_1 -2l_A\eps)a_\eps.
\end{align*}
We define $W^d_i = (0,l_1) \times W_i$ and the energy
$${\cal E}^i_\eps(y) := \eps^{d-1} \sum_{\bar{x} \in (\calL'_\eps(V^d_i))^\circ} W_{\rm cell}(\bar{\nabla} y(\bar{x})) .$$
For $C_1 \geq 2a_{\rm crit}$, $C_2 > 0$ sufficiently large and for $\delta> 0$ as before choose $f$ as in Theorem \ref{th: meso}. Then for $\eps$ small enough 
\begin{align*}
{\cal E}^i_\eps(y) & \ge \eps^{d-1}\sum^n_{j=1} \Big( \sum_{\bar{x} \in (\calL'_\eps(U_j))^\circ} W_{\rm cell}(\bar{\nabla} y(\bar{x})) + \frac{1}{2} \sum_{\bar{x} \in \partial_{W_i}(\calL'_\eps(U_j))} W_{\rm cell}(\bar{\nabla} y(\bar{x})) \Big) \\
&  \geq \sum^n_{j=1} M(U_j, r_j) \geq \sum^n_{j=1} f(r_j,u_j - u_{j-1}).
\end{align*}
Here we observe that due to the construction of the sets $W_i$ we have $\partial_{W_i}(\calL'_\eps(U_j)) \subset (\calL'_\eps(V^d_i))^\circ$ for all $j=1,\ldots,n$, $i \in I$. If there is some $j$ such that $r_j \geq C_2  \max \lbrace 1,  \eps^{(s - 1)(3 - d)}\rbrace$, then ${\cal E}^{i}_\eps(y) \geq \tilde{\omega}(|W_i|) \big(\frac{\beta_A}{\det A} -\delta\big)$ by \eqref{eq: sharp2}, where $\tilde{\omega}(\cdot) = \eps^{s(d-1)} \omega(\eps^{-s(d-1)} \cdot)$ (note that now $|W_i| \sim \eps^{s(d-1)}$). Otherwise all $r_j$ lie in the regime where $f$ is convex in $r$ and linear in $\lambda$. We then compute using Jensen's inequality 
\begin{align*}
  {\cal E}^{i}_\eps(y) 
  & \geq \sum^n_{j=1} f(r_j,u_j - u_{j-1}) = \sum^n_{j=1} \frac{u_j - u_{j-1}}{\lambda_0\eps^s} f(r_j,\lambda_0 \eps^s) \\
  & \geq \frac{\sum^n_{j=1} u_j - u_{j-1}}{\lambda_0 \eps^s} \  f\bigg(\frac{\sum^n_{j=1}  (u_j - u_{j-1}) \, r_j}{\sum^n_{j=1} u_j - u_{j-1}} ,\lambda_0 \eps^s\bigg)  \\
  & = \frac{l_1  - 2 l_A \eps}{\lambda_0 \eps^s} \ f\Big(a_\eps, \lambda_0\eps^s \Big), 
\end{align*}
whence for $a_\eps \geq 2a_{\rm crit} \sqrt{\eps}$, due to the monotonicity of $f$, also 
\begin{align*}
  {\cal E}^{i}_\eps(y) 
  &\geq \frac{l_1  - 2 l_A \eps}{\lambda_0 \eps^s} \ f\Big(2a_{\rm crit} \sqrt{\eps}, \lambda_0\eps^s \Big) \\ 
  &= \tilde{\omega}(|W_i|) \Big(\frac{(l_1  - 2 l_A \eps) \alpha_A 4 a_{\rm crit}^2}{2\det A} -   (l_1 - 2 l_A \eps) \delta\Big) 
  \geq \tilde{\omega}(|W_i|) \Big(\frac{\beta_A}{\det A} -\delta\Big)
\end{align*}
by \eqref{eq: sharp1}, where the last inequality holds for $\delta$ small enough. Repeating the calculation for $a_\eps \leq 2a_{\rm crit}\sqrt{\eps}$ and using \eqref{eq: sharp1} yields
$${\cal E}^{i}_\eps(y)  \geq \tilde{\omega}(|W_i|) m_\eps :=  \tilde{\omega}(|W_i|) \min\Big\lbrace \frac{(l_1  -2l_A\eps)  \alpha_A a_\eps^2}{2\eps\det A}  - (l_1  -2l_A\eps) \delta, \frac{\beta_A}{\det A} -\delta \Big\rbrace.$$ 
Using that ${\cal E}_\eps(y) \geq \sum_{i \in I} {\cal E}^i_\eps(y) $ and $\tilde{\omega}(|W_i|) \geq \sigma(\delta) \eps^{s(d-1)}$, where $\sigma(\delta)= \min\lbrace \omega(s): 1- C\delta^{\frac{d-1}{d-2}} \le s \le 1\rbrace \le 1$ for all $i \in I$ and recalling \eqref{eq: broken energy} we get for $\delta$ small enough 
\begin{align*}
\liminf_{\eps \to 0} \inf \lbrace {\cal E}_\eps(y): y\in {\cal A}(a_\eps)\rbrace & \geq \liminf_{\eps \to 0} \Big(\# \bar{I} \  \frac{\eps^{s(d-1)} 2 \beta_A}{\det A} + \# (I \setminus \bar{I}) \,\sigma(\delta) \eps^{s(d-1)} \ m_\eps\Big) \\ 
& \geq \liminf_{\eps \to 0} \# I  \ \sigma(\delta) \eps^{s(d-1)} \ m_\eps \\ 
& \geq \sigma(\delta) \,  \prod\limits^d_{j=2} l_j \min\Big\lbrace \frac{l_1 \alpha_A a^2}{2\det A} -  l_1 \delta, \frac{\beta_A}{\det A} -\delta \Big
\rbrace,
\end{align*}
as $a_\eps/\sqrt{\eps} \to a$. Letting $\delta \to 0$ shows
$$ \liminf_{\eps \to 0} \inf \lbrace {\cal E}_\eps(y): y\in {\cal A}(a_\eps)\rbrace  \geq  \frac{ \prod^d_{j=2} l_j}{\det A} \min \Big\lbrace \frac{1}{2} l_1 \alpha_A a^2, \beta_A \Big\rbrace.$$
Here we used that $\lim_{\delta \to 1} \sigma (\delta) = 1$. It remains to prove that the right hand side in Theorem \ref{th: cleavage} is attained for some sequence of deformations. This essentially follows from the sharpness of the estimates \eqref{eq: sharp1} and \eqref{eq: sharp2}. In particular, as in \eqref{eq: el min} for $a < \infty$ we consider
 \begin{align}\label{eq: elastic min}
y^{\rm el}_\eps(x) = x +  \bar{F}(a_\eps) \, x, \ \ \ x \in \calL_{\eps} \cap \Omega,
\end{align}
and as in the proof of Lemma \ref{lemma: el} it is not hard to see that
$$ \lim_{\eps \to 0} {\cal E}_\eps (y^{\rm el}_\eps) =  \prod\limits^d_{j=1} l_j \frac{\alpha_A}{2\det A} \lim_{\eps \to 0}\Big(\frac{a_\eps}{\eps}\Big)^2 =  \prod\limits^d_{j=1} l_j \frac{\alpha_A}{2\det A} a^2.$$
For $y^{\rm cr}_\eps$ we proceed as in \eqref{eq: cr min}: We choose $\xi$ such that \eqref{eq: beta def} is satisfied. Due to the assumption $l_1 \geq L$ it is possible to define a hyperplane $\Pi = \lbrace x \in \R^d: x \cdot \xi = c\rbrace$ such that $\Pi \cap \overline{\Omega} \subset \overline{\Omega}\setminus (B_1^\eps \cup B^\eps_2)$. We let
\begin{align}\label{eq: fracture min}
y^{\rm cr}_\eps(x) = \begin{cases}  x, &  x\cdot \xi < c, \\
                      x + l_1 a_\eps  \e_1,   & x\cdot \xi > c,    \end{cases} \ \ \ x \in \calL_{\eps} \cap \Omega. 
\end{align}
Again counting the quantity of broken springs as in \eqref{eq: counting} we derive $ \lim_{\eps \to 0} {\cal E}_\eps (y^{\rm cr}_\eps) =  \prod\limits^d_{j=2} l_j \frac{\beta_A }{\det A}$. \eop

\section{Examples: mass-spring models}\label{sec: examples}

In the following we examine several mass-spring models to which the above results apply. We calculate the constants $\alpha_A$, $\beta_A$ explicitly and thus we can provide the limiting energy of Theorem \ref{th: cleavage} as well as the critical value of boundary displacements $a_{\rm crit.}$ Moreover, we specify minimizing configurations and discuss their behavior depending on the properties of the cell energy. 

Note that the cell energies under consideration which consist of pair interaction energies are typically minimized on $\bar{O}(d)$. Thus, parts of the specimen might flip their orientation without affecting the energy. In order avoid such unphysical behavior and to satisfy Assumption \ref{assu: 1}(i) we introduce a frame indifferent penalty term $\chi \geq 0$ vanishing in a neighborhood of $\bar{SO}(d)$ and $\infty$ and satisfying $\chi \geq c_\chi > 0$ in a neighborhood of $\bar{O}(d) \setminus \bar{SO}(d)$. For example, in line with a widely used local orientation preserving condition in the continuum setting we may set 
$$ \chi(\bar{\nabla} y(\bar{x})) 
   = \begin{cases} 
       0, 
       & \text{if} \  \det(\nabla \tilde{y}) > 0 \text{ a.e.\ on } Q(\bar{x}) \text{ or} \ |\bar{\nabla}y(\bar{x})|\geq R \\ 
       \infty 
       &  \text{otherwise} 
     \end{cases} $$
for some $R \gg 1$. The penalty term does not change the energy in the elastic and fracture regime.

\subsection{Triangular lattices with NN interaction}

We begin with a planar model where the atoms in the reference configuration are given by the portion of a triangular lattice lying in $\Omega = (0,l_1) \times (0,l_2)$ and only interact with their nearest neighbors. It serves as the most basic non-trivial example to which our theory applies. In fact, a much more complete analysis of this model including a detailed characterization of low energy configurations has been performed in \cite{FriedrichSchmidt:2011}. Let $\calL = A\Z^2 = T_\phi \begin{pmatrix} 1 & \frac{1}{2} \\  0 &  \frac{\sqrt{3}}{2}\end{pmatrix} \Z^2$ with $T_\phi = \begin{pmatrix}  \cos \phi & - \sin\phi \\  \sin\phi & \cos\phi \end{pmatrix}$ for $\phi \in [0,\frac{\pi}{3})$ and for a deformation $y: \calL_\eps \cap \Omega \to \R^2$ let
$${\cal E}_\eps (y) = \frac{\eps}{2} \sum_{x,x' \in \calL_\eps  \atop |x - x'| = \eps} W \Big(\frac{|y(x) - y(x')|}{\eps} \Big) + \eps \sum_{\bar{x} \in (\calL'_\eps(\Omega))^\circ}\chi(\bar{\nabla} y(\bar{x})),$$
where $W:[0,\infty) \to [0,\infty)$ is of `Lennard-Jones-type', i.e.\ (1) $W \geq 0$, $W(r) = 0 \Leftrightarrow r=1$, (2) $W$ is continuous and $C^2$ in a neighborhood of $1$ with $\bar{\alpha} :=  W''(1) > 0$, (3) $\lim_{r \to \infty} W(r) = \bar{\beta} > 0$. Denoting the $i$-th column of $G$ by $G_i$ and letting $Z = A \, \frac{1}{2} \begin{pmatrix}  -1 & 1& 1&  -1 \\ -1 & -1& 1&  1   \end{pmatrix}$ the cell energy can be written as 
\begin{align*}
W_{\rm cell}(G) & = \frac{1}{2} \big( W(|G_2 - G_1|) + W(|G_3 - G_2|) +W(|G_4 - G_3|)  \\ &\quad +W(|G_1 - G_4|)  +2W(|G_4 - G_2|) + \chi(G) \big).
\end{align*}
In \cite{FriedrichSchmidt:2011} is shown that this is an admissible cell energy in the sense of Assumption \ref{assu: 1}. We compute
$$ {\cal Q} = \frac{3\bar{\alpha}}{8}\begin{pmatrix} 3 & 1 & 0 \\  1 & 3 & 0  \\ 0 & 0 & 2 \end{pmatrix}$$
and therefore $\alpha_A = \bar{\alpha}$. The fact that ${\cal Q}$ is independent of $\phi$ particularly shows that the energy in the linearized elastic regime is isotropic. With $\nu^\phi_1 = T_\phi(1, 0)^T$, $\nu^\phi_2 = T_\phi(\frac{1}{2},\frac{\sqrt{3}}{2})^T$ and $\nu^\phi_3 = T_\phi(-\frac{1}{2},\frac{\sqrt{3}}{2})^T$ we get
$$\beta_A = \min_{\varsigma \in S^1} \bar{\beta} \frac{\sum^3_{i=1} |\nu^\phi_i \cdot \varsigma|}{|\e_1 \cdot \varsigma|} = \frac{\sqrt{3} \bar{\beta}}{ \sin(\phi+\frac{\pi}{3})}$$
and then we re-derive
$$ {\cal E}_{\rm lim}(a) 
   = \frac{2 l_2}{\sqrt{3}} \min \Big\lbrace \frac{1}{2} l_1 \bar{\alpha} a^2, \frac{\sqrt{3}\bar{\beta}}{ \sin(\phi+\frac{\pi}{3})}\Big\rbrace, \ \ \ a_{\rm crit} 
   = \sqrt{\frac{2 \sqrt{3} \bar{\beta}}{l_1\bar{\alpha} \, \sin(\phi+\frac{\pi}{3})}}.$$
The minimizers in the elastic regime have the form \eqref{eq: elastic min} with 
$$\bar{F}(a_\eps) = \begin{pmatrix}
 a_\eps & 0 \\ 0 & - \frac{a_\eps}{3}\end{pmatrix}$$
and in particular show the Poisson effect. In the supercritical case deformations being cleaved along lines in $\nu^\phi_2$-direction are energetically optimal. In fact, it can be shown that except for the symmetrically degenerate case $\phi = 0$, for which both cleavage in $\nu^\phi_2$-direction and $\nu^\phi_3$-direction is optimal, the rescaled displacements of deformations with almost optimal energy converge, up to suitable translations and shifts in $c$, to $(a l_1, 0) \chi_{\{x \in \Omega : x \cdot \xi > c\}}$, see \cite{FriedrichSchmidt:2011}.

\subsection{Square lattices with NN and NNN interaction}

The behavior in the elastic regime of the following two dimensional model comprising nearest and next to nearest neighbor atomic interactions was treated by Friesecke and Theil in \cite{FrieseckeTheil:02}. We let $\calL = A \Z^2 = T_\phi \Z^2$ for $\phi \in [0,\frac{\pi}{2})$ and
\begin{align*}
{\cal E}_\eps(y) 
&=  \frac{\eps}{2} \sum_{x,x' \in \calL_\eps  \atop |x - x'| = \eps} W_1 \Big(\frac{|y(x) - y(x')|}{\eps} \Big) \\ 
&+ \frac{\eps}{2} \sum_{x,x' \in \calL_\eps  \atop |x - x'| = \sqrt{2}\eps} W_2 \Big(\frac{|y(x) - y(x')|}{\sqrt{2}\eps} \Big)  + \eps \sum_{\bar{x} \in (\calL'_\eps(\Omega))^\circ}\chi(\bar{\nabla} y(\bar{x})) 
\end{align*}
for deformations $y: \calL_\eps \cap \Omega \to \R^2$ and potentials $W_1$, $W_2$ as above with $\bar{\alpha}_1$, $\bar{\beta}_1$ and $\bar{\alpha}_2$, $\bar{\beta}_2$, respectively. The associated cell energy is given by
\begin{align*}
W_{\rm cell}(G) = \frac{1}{4} \sum_{|z_i - z_j| = 1} W_1(|G_i - G_j|) + \frac{1}{2} \sum_{|z_i - z_j| = \sqrt{2}} W_2\Big(\frac{|G_i - G_j|}{\sqrt{2}}\Big) + \chi(G)
\end{align*} 
In \cite{FrieseckeTheil:02} it is shown that the cell energy is admissible. We calculate
$${\cal Q} = \frac{1}{2}T^{*T}_\phi\begin{pmatrix}  2\bar{\alpha}_1  + \bar{\alpha}_2 & \bar{\alpha}_2 & 0 \\  \bar{\alpha}_2 & 2 \bar{\alpha}_1  + \bar{\alpha}_2 & 0  \\ 0 & 0 & 2\bar{\alpha}_2 \end{pmatrix} T^*_\phi,$$
with
$$T^*_\phi = \begin{pmatrix} c^2 & s^2 & -\sqrt{2}cs\\  s^2 & c^2 & \sqrt{2}cs  \\ \sqrt{2}cs & -\sqrt{2}cs & c^2 -s^2 \end{pmatrix},  \ \ \ c := \cos\phi, \ s := \sin\phi.$$
An elementary computation then shows
$$\alpha_A = \frac{2\bar{\alpha}_1 \bar{\alpha}_2 (\bar{\alpha}_1 +  \bar{\alpha}_2) }{2\bar{\alpha}_1 \bar{\alpha}_2 + 4 \bar{\alpha}^2_1 c^2 s^2 + \bar{\alpha}^2_2 (c^2 - s^2)^2}.$$
Letting $\nu^\phi_1 = T_\phi \e_1$, $\nu^\phi_2 = T_\phi \e_2$, $\nu^\phi_3 = T_\phi (\e_1 + \e_2)$ and $\nu^\phi_4 = T_\phi (\e_1 - \e_2)$ and $\gamma_1 = \max\lbrace c, s\rbrace$, $\gamma_2 = c + s$ we obtain for the fracture constant 
\begin{align*}
\beta_A  & = \min_{\varsigma \in S^1} \frac{\bar{\beta_1} ( |\nu^\phi_1 \cdot \varsigma| +  |\nu^\phi_2 \cdot \varsigma| )  +  \bar{\beta_2} ( |\nu^\phi_3 \cdot \varsigma| +  |\nu^\phi_4 \cdot \varsigma| ) }{|\e_1 \cdot \varsigma|} \\
&  = \min\Big\lbrace\frac{\bar{\beta}_1 + 2\bar{\beta}_2}{\gamma_1}, \frac{ 2 \bar{\beta}_1 + 2 \bar{\beta}_2}{\gamma_2}\Big\rbrace.
\end{align*}
Here we used that it suffices to minimize over the set ${\cal P} = \lbrace \frac{\varsigma}{|\varsigma|}: \varsigma = \nu_i^\phi, i=1,\ldots,4 \rbrace \subset S^1$. 
Below the critical value $a_{\rm crit}$ energetically optimal configurations are given by functions of the form \eqref{eq: elastic min} with 
$$\bar{F}(a_\eps) =  \begin{pmatrix}
a_\eps  & 0  \\        \frac{(\bar{\alpha}^2_1 - \bar{\alpha}^2_1) c s (c^2-s^2) }{2\bar{\alpha}_1\bar{\alpha}_2 + \bar{\alpha}^2_2 (c^2 -s^2)^2 + 4\bar{\alpha}^2_1 c^2 s^2}a_\eps               & \frac{-\bar{\alpha}^2_2 (c^2-s^2)^2 - 4\bar{\alpha}^2_1 c^2s^2}{2\bar{\alpha}_1 \bar{\alpha}_2 + \bar{\alpha}^2_2 (c^2-s^2)^2 + 4\bar{\alpha}^2_1c^2s^2} a_\eps     
\end{pmatrix}. $$
In particular, the configurations show the Poisson effect and in the case that $\bar{\alpha}_1 \neq \bar{\alpha}_2$ and $\phi \in (0, \frac{\pi}{2}) \setminus \lbrace \frac{\pi}{4}\rbrace$ also shear effects occur. Limiting minimal configurations beyond critical loading are given by deformations of the form \eqref{eq: fracture min}, where the normal $\xi$ to the hyperplane $\Pi$ is an element of $\lbrace\nu^\phi_i: i=1,\ldots 4\rbrace$. While in the previous example the cleavage direction was determined only by the geometry of the problem (i.e. by $\phi$), it now depends also on the ratio of $\bar{\beta}_1$, $\bar{\beta}_2$.

We note that here for every $\phi \in (0, \frac{\pi}{2}) \setminus \lbrace \frac{\pi}{4}\rbrace$ by choosing the specific values $\bar{\beta}_1 = 1$ and $\bar{\beta}_2 = \frac{1}{2}\max\{\cot \phi, \tan \phi\} - \frac{1}{2}$ the minimum in the expression for $\beta_A$ is attained at $\nu^\phi_2$ and $\nu^\phi_3$, respectively, at $\nu^\phi_1$ and $\nu^\phi_4$. As a consequence, unlike for the triangular lattice in the previous example, also for general lattice orientations there may be deformations with almost optimal energy whose rescaled displacements in the continuum limit have a serrated jump set.

\subsection{Cubic lattices with NN and NNN interaction}

We consider the following three dimensional model with nearest and next nearest interactions in the domain $\Omega = (0,l_1) \times (0,l_2) \times (0,l_3)$. We let $\calL = A \Z^3 = T_{\phi,\psi} \Z^3$, where 
$$T_{\phi,\psi} = \begin{pmatrix} \cos\psi&  -\sin\psi & 0\\ \sin\psi & \cos\psi & 0 \\ 0 & 0 & 1 \end{pmatrix}\begin{pmatrix}  1 & 0 & 0 \\ 0 & \cos\phi & -\sin\phi  \\ 0 & \sin\phi & \cos\phi  \end{pmatrix}  \ \text{for} \ \phi,\psi \in [0,\frac{\pi}{2}).$$ 
We let
\begin{align*}
{\cal E}_\eps(y) &=  \frac{\eps^2}{2} \sum_{x,x' \in \calL_\eps  \atop |x - x'| = \eps} W_1 \Big(\frac{|y(x) - y(x')|}{\eps} \Big) \\ 
&+ \frac{\eps^2}{2} \sum_{x,x' \in \calL_\eps  \atop |x - x'| = \sqrt{2}\eps} W_2 \Big(\frac{|y(x) - y(x')|}{\sqrt{2}\eps} \Big)  +  \eps^2 \sum_{\bar{x} \in (\calL'_\eps(\Omega))^\circ}\chi(\bar{\nabla} y(\bar{x})) 
\end{align*}
for deformations $y: \calL_\eps \cap \Omega \to \R^2$ and potentials $W_1$, $W_2$ as above with $\bar{\alpha}_1$, $\bar{\beta}_1$ and $\bar{\alpha}_2$, $\bar{\beta}_2$, respectively. The associated cell energy is given by
\begin{align*}
W_{\rm cell}(G) = \frac{1}{8} \sum_{|z_i - z_j| = 1} W_1(|G_i - G_j|) + \frac{1}{4} \sum_{|z_i - z_j| = \sqrt{2}} W_2\Big(\frac{|G_i - G_j|}{\sqrt{2}}\Big) + \chi(G).
\end{align*} 
In \cite{Schmidt:2006} it has been shown that the cell energy is admissible. As before, an elementary computation shows 
$${\cal Q} =  T^{*T}_{\psi} T^{*T}_{\phi}\frac{1}{2}\begin{pmatrix}
 2\bar{\alpha}_1 + 2\bar{\alpha}_2 &  \bar{\alpha}_2   &  \bar{\alpha}_2  & 0 & 0  &  0 \\
\bar{\alpha}_2 & 2\bar{\alpha}_1 + 2\bar{\alpha}_2   &  \bar{\alpha}_2  & 0 & 0  &  0  \\
\bar{\alpha}_2 &   \bar{\alpha}_2 & 2\bar{\alpha}_1 + 2\bar{\alpha}_2 & 0 & 0  &  0   \\ 
0 & 0 & 0 & 2 \bar{\alpha}_2 & 0 & 0 \\
0 & 0 & 0 & 0 &2 \bar{\alpha}_2  & 0 \\
0 & 0 & 0 & 0 & 0 & 2 \bar{\alpha}_2  
\end{pmatrix} T^*_{\phi} T^*_{\psi} ,$$
where 
$$T^*_\phi =
\begin{pmatrix}
1 & 0 & 0 & 0 & 0 & 0 \\
0 & c_1^2 & s_1^2 & 0 & 0 & -2c_1s_1 \\
0 & s_1^2 & c_1^2 & 0 & 0 & 2c_1s_1 \\
0 & 0 & 0 & c_1 & -s_1 & 0 \\
0 & 0 & 0 & s_1 & c_1 & 0 \\
0 & c_1s_1 & -c_1s_1 & 0 & 0 &c_1^2 - s_1^2
\end{pmatrix},
T^*_\psi =
\begin{pmatrix}
c_2^2 &  s_2^2 & 0 & -2 c_2 s_2 & 0 & 0  \\
s_2^2 &  c_2^2 & 0 & 2 c_2 s_2 & 0 & 0  \\
0 & 0 & 1 & 0 & 0 & 0 \\
c_2 s_2 & -c_2 s_2 & 0 & c_2^2 - s_2^2 & 0 & 0 \\
0 & 0 & 0 & 0 & c_2 & -s_2 \\
0 & 0 & 0 & 0 & s_2 & c_2
\end{pmatrix},
$$
with the abbreviations $c_1 = \cos\phi$, $c_2 = \cos\psi$, $s_1 = \sin\phi$ and $s_2 = \sin\psi$. Applying Lemma \ref{lemma: alpha Q} we then obtain
$$\alpha_A = \frac{\bar{\alpha}_2 (2\bar{\alpha}_1 + \bar{\alpha}_2)^2(\bar{\alpha}_1 + 2\bar{\alpha}_2 )}{
8\bar{\alpha}_1^3 c_2^2 s_2^2 + 2\bar{\alpha}_1\bar{\alpha}_2^2(4 - c_2^2 s_2^2) + 
  4\bar{\alpha}_1^2\bar{\alpha}_2(4 c_2^2 s_2^2 - 1) + \bar{\alpha}_2^3(3 - 4c_2^2 s_2^2)}. $$
In particular, $\alpha_A$ is independent of $c_1$ and $s_1$. We let
${\cal V}_1^{\phi,\psi} = T_{\phi,\psi} \lbrace  \e_1,\e_2,\e_3 \rbrace$, ${\cal V}_2^{\phi,\psi} = T_{\phi,\psi} \lbrace  \e_1 + \e_2 , \e_1 - \e_2 , \e_1 + \e_3, \e_1 - \e_3, \e_2 + \e_3, \e_2 - \e_3\rbrace$  and $\gamma_1 = \max\lbrace |c_2|, |c_1s_2|, |s_1s_2|\rbrace$, $\gamma_2 = \max \lbrace |c_2 \pm c_1s_2|, |c_2 \pm s_1s_2|, |c_1s_2 \pm s_1s_2|\rbrace$, $\gamma_3 = \max\lbrace |c_2 \pm c_1s_2 \pm s_1s_2|\rbrace$, $\gamma_4 = \max\lbrace |2c_2 \pm c_1s_2 \pm s_1s_2|, |c_2 \pm 2c_1s_2 \pm s_1s_2|, |c_2 \pm c_1s_2 \pm 2s_1s_2|\rbrace$. One can show that 
\begin{align*}
{\cal P} &= \lbrace \varsigma/|\varsigma| : \varsigma = T_{\phi,\psi}\e_1, i=1,2,3 \rbrace  \cup \lbrace \varsigma/|\varsigma| : \varsigma = T_{\phi,\psi} (\e_1 \pm \e_2 \pm \e_3) \rbrace \\ 
&  \ \ \ \cup \lbrace \varsigma/|\varsigma| : \varsigma = T_{\phi,\psi} (\e_1 \pm \e_2 \pm \e_3 \pm \e_i), i=1,2,3 \rbrace.
\end{align*} 
Then
\begin{align*}
\beta_A &= \min_{\varsigma \in S^1} \frac{\sum_{\nu \in {\cal V}_1^{\phi,\psi}} \bar{\beta}_1 |\nu \cdot \varsigma| + \sum_{\nu \in {\cal V}_2^{\phi,\psi}} \bar{\beta}_2 |\nu \cdot \varsigma|}{|\e_1 \cdot \varsigma|} \\
& = \min \Big\lbrace\frac{\bar{\beta}_1 + 4\bar{\beta}_2}{\gamma_1},\frac{2\bar{\beta}_1 + 6\bar{\beta}_2}{\gamma_2},\frac{3\bar{\beta}_1 + 6\bar{\beta}_2}{\gamma_3}, \frac{4\bar{\beta}_1 + 10\bar{\beta}_2}{\gamma_4}\Big\rbrace.
\end{align*}


 \typeout{References}

\end{document}